\documentclass{conm-p-l}

\usepackage[dvips]{graphics}
\usepackage{ifthen}
\usepackage{soul}
\usepackage{latexsym}
\usepackage{amsthm,amssymb}
\usepackage{amsbsy,amsfonts,amsmath}

\usepackage{a4}

\numberwithin{equation}{section}

\newtheorem{theorem}{Theorem}[section]

\newtheorem{lemma}[theorem]{Lemma}

\theoremstyle{definition}
\newtheorem{definition}[theorem]{Definition}

\theoremstyle{remark}
\newcommand{\al}{\alpha}

\newcommand{\tq}{q}

\newcommand{\defdimofV}{\theta}
\newcommand{\bQ}{{\mathbb{Q}}}
\newcommand{\bR}{{\mathbb{R}}}
\newcommand{\bZ}{{\mathbb{Z}}}
\newcommand{\bN}{{\mathbb{N}}}
\newcommand{\bJ}{{\mathbb{J}}}
\newcommand{\bO}{{\mathbb{O}}}
\newcommand{\bZp}{\bZ_{\geq 0}}

\newcommand{\bA}{{\mathbb{A}}}

\newcommand{\defI}{I}

\newcommand{\bK}{{\mathbb{K}}}
\newcommand{\bKt}{{\bK^\times}}

\newcommand{\defV}{V}
\newcommand{\defVZ}{\defV_\bZ}
\newcommand{\bhm}{\chi}

\newcommand{\defPi}{\pi}%{\Pi}
\newcommand{\defU}{U}
\newcommand{\defUbhmPi}{{\defU(\bhm,\defPi)}}%{{\defU(\bhm,\Pi)}}
\newcommand{\defUo}{\defU^0}
\newcommand{\defUp}{\defU^+}
\newcommand{\defUm}{\defU^-}

\newcommand{\defUbhmopPi}{{\defU(\bhmop,\defPi)}}

\newcommand{\defacK}{K}
\newcommand{\defacL}{L}
\newcommand{\defacE}{E}
\newcommand{\defacF}{F}

\newcommand{\defacEc}{\defacE^\vee}
\newcommand{\defacFc}{\defacF^\vee}

\newcommand{\BdefacE}{{\bar{\defacE}}}
\newcommand{\BdefacF}{{\bar{\defacF}}}
\newcommand{\BdefacH}{{\bar{H}}}

\newcommand{\defUO}{\defU_\bO}
\newcommand{\defUoO}{\defUo_\bO}
\newcommand{\defUpO}{\defUp_\bO}
\newcommand{\defUmO}{\defUm_\bO}

\newcommand{\gravedefU}{{\grave{\defU}}}
\newcommand{\gravedefacK}{{\grave{\defacK}}}
\newcommand{\gravedefacL}{{\grave{\defacL}}}
\newcommand{\gravedefacE}{{\grave{\defacE}}}
\newcommand{\gravedefacF}{{\grave{\defacF}}}

\newcommand{\defVZPip}{\defVZ^{\defPi,+}}

\newcommand{\sqrtbhm}{\sqrt{\bhm}}

\newcommand{\dotq}{{\dot{q}}}
\newcommand{\ddotq}{{\ddot{q}}}

\newcommand{\defOmega}{\Omega}
\newcommand{\OmegabhmPi}{\defOmega^{\bhm,\defPi}}

\newcommand{\defGamma}{\Gamma}
\newcommand{\GammabhmPi}{\defGamma^{\bhm,\defPi}}
\newcommand{\GammabhmtauiPi}{\defGamma^{\bhm,\tauidefPi}}
\newcommand{\GammabhmopPi}{\defGamma^{\bhmop,\defPi}}

\newcommand{\bhmop}{\bhm^{\mathrm{op}}}
\newcommand{\defUpsilon}{\Upsilon}
\newcommand{\defUpbhmopPi}{\defUpsilon^{\bhmop,\defPi}}

\newcommand{\defRbhm}{R_\bhm}
\newcommand{\defRbhmpip}{\defRbhm^{\defPi,+}}
\newcommand{\defRbhmpi}{\defRbhm^\defPi}

\newcommand{\defvphibhmpip}{\varphi_\bhm^{\defPi,+}}
\newcommand{\defvphibhmpi}{\varphi_\bhm^\defPi}

\newcommand{\Nbhmpi}{N^{\bhm,\defPi}}

\newcommand{\deftaubhm}{\tau^\bhm}

\newcommand{\tauidefPi}{{\deftaubhm_i\defPi}}

\newcommand{\defvphibhmtauipi}{\varphi_\bhm^\tauidefPi}

\newcommand{\NbhmdtiPi}{N^{\bhm,\tauidefPi}}

\newcommand{\defRtaubhmipip}{\defRbhm^{\tauidefPi,+}}
\newcommand{\defRtaubhmipi}{\defRbhm^\tauidefPi}

\newcommand{\defUbhmtauiPi}{{\defU(\bhm,\tauidefPi)}}

\newcommand{\LusT}{T}

\newcommand{\smallT}{\xi}

\newcommand{\rmad}{{\mathrm{ad}}}

\newcommand{\defTheta}{\Theta}

\newcommand{\spa}{\varpi}

\newcommand{\HopfD}{\Delta}
\newcommand{\HopfS}{S}
\newcommand{\Hopfe}{\varepsilon}
\newcommand{\tHopfD}{\HopfD}
\newcommand{\tHopfS}{\HopfS}
\newcommand{\tHopfe}{\Hopfe}
\newcommand{\tBp}{\defU^{+,\flat}}
\newcommand{\tBm}{\defU^{-,\flat}}
\newcommand{\tBpbhmPi}{\tBp(\bhm,\defPi)}
\newcommand{\tBmbhmPi}{\tBm(\bhm,\defPi)}
\newcommand{\rmSpan}{{\mathrm{Span}}}

\newcommand{\tvt}{\vartheta}
\newcommand{\tvtbhmPi}{\vartheta^{\bhm,\defPi}}
\newcommand{\rmid}{{\mathrm{id}}}

\newcommand{\mbhmpi}{m^{\bhm,\defPi}}

\newcommand{\LusTi}{\LusT^{\bhm,\tauidefPi}_i}
\newcommand{\LusTiop}{\LusT^{\bhmop,\tauidefPi}_i}
\newcommand{\LusTidefPi}{\LusT^{\bhm,\defPi}_i}

\newcommand{\HdefacE}{{\hat{\defacE}}}

\newcommand{\rmidbhmPiT}{\rmid^{\bhm,\defPi}\LusT}

\newcommand{\zetabhmPi}{\zeta^{\bhm,\defPi}}

\newcommand{\defQp}{Q_1}
\newcommand{\defQm}{Q_2}
\newcommand{\defQpm}{Q_t}

\newcommand{\mfkgp}{{\mathfrak{g}}_1}
\newcommand{\mfkgm}{{\mathfrak{g}}_2}
\newcommand{\mfkgpm}{{\mathfrak{g}}_t}

\newcommand{\barmfkgp}{{\bar {\mathfrak{g}}}^+}
\newcommand{\barmfkgm}{{\bar {\mathfrak{g}}}^-}

\newcommand{\defUA}{\defU_\bA}
\newcommand{\defUpA}{\defUp_\bA}
\newcommand{\defUmA}{\defUm_\bA}
\newcommand{\defUoA}{\defUo_\bA}

\newcommand{\tW}{W}
\newcommand{\ts}{s}
\newcommand{\te}{e}
\newcommand{\defell}{\ell}

\newcommand{\twoh}{w_\circ}

\newcommand{\ellwo}{\defell(\twoh)}

\begin{document}

\title{Kostant-Lusztig $\bA$-bases of Multiparameter Quantum Groups}

\author[N. Jing]{Naihuan Jing}
\address{Department of Mathematics, North Carolina State University,  Raleigh,
NC 27695-8205, USA}
\email{jing@ncsu.edu}
\thanks{NJ is partially supported by National Natural Science Foundation of China grant \# 11531004 and Simons Foundation grant \# 523868.}

\author[K.C. Misra]{ Kailash C. Misra}
\address{Department of Mathematics, North Carolina State University,  Raleigh,
NC 27695-8205, USA}
\email{misra@ncsu.edu}
\thanks{KCM is partially supported by Simons Foundation grant \# 307555.}

\author[H. Yamane]{Hiroyuki Yamane}
\address{Department of Mathematics, Faculty of Science, University of Toyama, Gofuku, Toyama 930-8555, Japan}
\email{hiroyuki@sci.u-toyama.ac.jp}
\thanks{HY is partially supported by JSPS Grand-in-Aid for Scientific Research (C), 16K05095.}

\subjclass[2010]{Primary 17B37,17B10; Secondary 81R50}

%\begin{center}
%{\Huge{Kostant-Lusztig $\bA$-base of Multiparameter Quantum Groups}}
%\end{center}
%
%\vspace{1cm}
%
%\begin{center}
%Naihuan Jing,\,\,Kailash C. Misra\,\,and\,\,Hiroyuki Yamane$^*$
%\end{center}

\begin{abstract}
We study the Kostant-Lusztig $\bA$-base of the multiparameter quantum groups.
To simplify calculations, especially for $G_2$-type,
we utilize the duality of the pairing of the universal $R$-matrix.
\end{abstract}

\maketitle

\section{Introduction}

Quantum enveloping algebras $U_q(\mathfrak g)$ \cite{Dr86, J} of the Kac-Moody algebras $\mathfrak g$ are one of the important classes of quantum groups. Quantum enveloping algebras and their integrable highest weight representations
enjoy favorable properties, among which the canonical bases of Lusztig \cite{Lus90b} or global crystal bases of Kashiwara \cite{K91} are the most prominent ones.
Existence of such canonical bases has also been established for quantum enveloping algebras $U_q(\mathfrak g)$  of Borcherds' generalized Kac-Moody algebras $\mathfrak g$ \cite{KKK}.

On the other hand, quantum enveloping algebras have been extended to Nichols algebras of diagonal types, which include multiparameter quantum enveloping algebras as examples, in particular, one-parameter quantum enveloping algebras $U_q(\mathfrak g)$ in Kac-Moody types. In \cite{AS} Andruskiuwitch and Schneider proved that finite dimensional pointed Hopf algebras with finite abelian group (with order $>7$) of group-like elements are essentially Lusztig's small quantum groups and their variants. Furthermore, Heckenberger classified the Nicholas algebras with arithmetic root data \cite{Hec09} and proved results similar to quantum enveloping algebras (see also \cite{Hec06}).

Lusztig~\cite{Lus90a} introduced the Kostant-Lusztig $\bA$-base of the quantum groups.
The $\bA$-base and the finite dimensional Hopf algbebra
of the quantum group at root of unity have been key figures in study of Lusztig conjectures, (see \cite{Jan03} for history).
In this paper we establish the Kostant-Lusztig $\mathbb A$-forms for the multiparameter quantum groups and construct invariant bases for each factor of the triangular decomposition.
Our general result is based on the structure theory for multiparameter quantum groups and information on the lower rank cases, most notably
the case of $G_2$. The recent work of Fan and Li \cite{FL} on two-parameter quantum algebras has made us believe that one should be able to use our $\mathbb A$-forms to construct canonical bases for the multiparameter quantum groups.

\section{Generalized quantum groups}
In this section we recall definitions and known results about multiparameter quantum groups.
We will use the following notations throughout this paper. The ring of real numbers and integers will be denoted by $\bR$ and $\bZ$ respectively. For $x$, $y\in\bR$, let $\bJ_{x,y}:=\{z\in\bZ|x\leq z\leq y\}$ and $\bJ_{x,\infty}:=\{z\in\bZ|x\leq z\}$.
Then $\bN=\bJ_{1,\infty}$ and $\bZp:=\bJ_{0,\infty}$. Let $\bK$ be a field of characteristic zero and
$\bKt:=\bK\setminus\{0\}$. For $x$, $y\in\bK$ and $r\in\bZp$, we denote $(r)_x:=\sum_{k=0}^{r-1}x^k$, $(r)_x!:=\prod_{k=1}^r(k)_x$, $(r;x,y):=1-x^{r-1}y$ and $(r;x,y)!:=\prod_{k=1}^r(k;x,y)$.

\subsection{Definition of generalized quantum groups}
Let $\defdimofV\in\bN$ and $I:=\bJ_{1,\defdimofV}$.
Let $\defV$ be a $\defdimofV$-dimensional $\bR$-linear space
with a basis $\{v_i|i\in\defI\}$.
Let $\defVZ:=\oplus_{i\in\defI}\bZ v_i$, so $\defVZ$ is a rank-$\defdimofV$ free $\bZ$-module
(or a free abelian group).
Let $\bhm:\defVZ\times\defVZ\to\bKt$ be a map
such that $\bhm(x+y,z)=\bhm(x,z)\bhm(y,z)$ and $\bhm(x,y+z)=\bhm(x,y)\bhm(x,z)$
hold for all $x$, $y$, $z\in\defVZ$.
We call such $\bhm$ {\it{a bi-character}}.

\begin{definition}\label{definition:defofGQG}
Let $\defPi:\defI\to\defVZ$ be a map such that
$\defPi(\defI)$ is a $\bZ$-base of $\defVZ$.
The generalized quantum group $\defU=\defUbhmPi$ is the unique associative $\bK$-algebra (with $1$) satisfying the following conditions (i)-(v).
\newline\newline
(i) As a $\bK$-algebra, $\defUbhmPi$ is generated by
$\defacK_\lambda$, $\defacL_\lambda$, ($\lambda\in\defVZ$), $\defacE_i$, $\defacF_i$
($i\in\defI$).
\newline\newline
(ii) The following equations hold. Let $\tq_{ij}:=\bhm(\defPi(i),\defPi(j))$.
\begin{equation*}
\begin{array}{l}
\defacK_0=1,\,\,\defacK_\lambda\defacK_\mu=\defacK_{\lambda+\mu}, \\
\defacL_0=1,\,\,\defacL_\lambda\defacL_\mu=\defacL_{\lambda+\mu},  \,\,
\defacK_\lambda\defacL_\mu=\defacL_\mu\defacK_\lambda,  \\
\defacK_{\defPi(i)}\defacE_j\defacK_{\defPi(i)}^{-1}=\tq_{ij}\defacE_j,\,\,\defacK_{\defPi(i)}\defacF_j\defacK_{\defPi(i)}^{-1}=\tq_{ij}^{-1}\defacF_j, \\
\defacL_{\defPi(i)}\defacE_j\defacL_{\defPi(i)}^{-1}=\tq_{ji}^{-1}\defacE_j,\,\,\defacL_{\defPi(i)}\defacF_j\defacL_{\defPi(i)}^{-1}=\tq_{ji}\defacF_j,  \\
\mbox{$[\defacE_i,\defacF_j]=\delta_{ij}(-\defacK_{\defPi(i)}+\defacL_{\defPi(i)})$}
\end{array}
\end{equation*}
(iii) There are subspaces $\defU_\lambda$ ($\lambda\in\defVZ$) such that
$\defU=\oplus_{\lambda\in\defVZ}\defU_\lambda$,
$\defU_\lambda\defU_\mu\subset\defU_{\lambda+\mu}$,
$\defacK_{\defPi(i)}^{\pm 1}\in\defU_0$, $\defacL_{\defPi(i)}^{\pm 1}\in\defU_0$,
$\defacE_i\in\defU_{\defPi(i)}$, $\defacF_i\in\defU_{-\defPi(i)}$
\newline\newline
(iv) Let $\defUo$ be the $\bK$-subalgebra of $\defU$ generated by
$\defacK_\lambda\defacL_\mu$ ($\lambda$, $\mu\in\defVZ$).
Let $\defUp$ (resp. $\defUm$) be the $\bK$-subalgebra of $\defU$ generated by
$\defacE_i$ (resp. $\defacF_i$)  ($i\in\defI$).
Then the elements $\defacK_\lambda\defacL_\mu$ ($\lambda$, $\mu\in\defVZ$)
form a $\bK$-basis of $\defUo$ and
we have the $\bK$-linear isomorphism
${\mathbf{m}}:\defUp\otimes\defUo\otimes\defUm\to\defU$ defined by
$X\otimes Y\otimes Z\mapsto XYZ$.
\newline\newline
(v) There exists no $X\in\defUp\setminus\{\mathbb{K}\}$
(resp. $Y\in\defUm\setminus\{\mathbb{K}\}$) such that
$[X,\defacF_i]=0$ ($[\defacE_i,Y]=0$) for all $i\in\defI$.
\end{definition}

Note that $\defUo\subset\defU_0$. For $\lambda\in\defVZ$,
let
$\defUp_\lambda:=\defUp\cap\defU_\lambda$ and $\defUm_\lambda:=\defUm\cap\defU_\lambda$.
Let $\defVZPip:=\oplus_{i\in I}\bZp\defPi(i)$.
Then $\defUp=\oplus_{\lambda\in\defVZPip}\defUp_\lambda$ and
$\defUm=\oplus_{\lambda\in\defVZPip}\defUm_{-\lambda}$.

\begin{lemma}\label{lemma:univ}
Let $\gravedefU$ be a $\bK$-algebra {\rm{(}}with $1${\rm{)}} generated by $\gravedefacK_\lambda$, $\gravedefacL_\lambda$, $\gravedefacE_i$
and $\gravedefacF_i$  satisfying conditions {\rm{(i)-(iv)}} above (in Definition~{\rm{\ref{definition:defofGQG}}}).
Then there exists a $\bK$-algebra epimorphism $\xi:{\grave{U}}\to\defU$ such that
$\xi(\gravedefacK_\lambda)=\defacK_\lambda$, $\xi(\gravedefacL_\lambda)=\defacL_\lambda$
$(\lambda\in\defVZ)$,
$\xi(\gravedefacE_i)=\defacE_i$, $\xi(\gravedefacF_i)=\defacF_i$
$(i\in\defI)$.
\end{lemma}
{\it{Proof.}} Let $\lambda\in\defVZPip$.
Assume that there exists $X\in\gravedefU^+_\lambda\setminus\{0\}$
such that $[X, \gravedefacF_i]=0$ for all $i\in\defI$.
Let ${\mathcal{I}}$ (resp. ${\mathcal{I}}^+$) be the two-sided ideal of $\gravedefU$
(resp. $\gravedefU^+$) generated by $X$. Then
${\mathcal{I}}={\rm{Span}}_\bK(\gravedefU^-\gravedefU^0{\mathcal{I}}^+)$.
Let $\gravedefU^\prime$ be the quotient $\bK$-algebra $\gravedefU/{\mathcal{I}}$.
Then $\gravedefU^\prime$ also satisfies the same conditions as
Definition~{\rm{\ref{definition:defofGQG}}}~{\rm{(i)-(iv)}}. Let $g:\gravedefU\to \gravedefU^\prime$ be the canonical map.
Then $g_{|\gravedefU^-}:\gravedefU^-\to(\gravedefU^\prime)^-$
and $g_{|\gravedefU^0}:\gravedefU^-\to(\gravedefU^\prime)^0$ are the $\bK$-algebra isomorphisms.
We see that $\dim (\gravedefU^\prime)^+_{\mu_1}= \dim \gravedefU^+_{\mu_1}$
for $\mu_1\in\defVZPip$ with $\mu_1-\lambda\notin\defVZPip$,
that $\dim (\gravedefU^\prime)^+_\lambda= \dim \gravedefU^+_\lambda-1$,
and that $\dim (\gravedefU^\prime)^+_{\mu_2}\leq\dim \gravedefU^+_{\mu_2}$
for $\mu_2\in\defVZPip$ with $\mu_2-\lambda\in\defVZPip$.
We also have a similar property for $\gravedefU^-$.
Then we can see the claim of this theorem by a standard argument using a direct limit. \par
\hfill $\Box$
\newline\par

By Lemma~\ref{lemma:univ}, we have the $\bK$-algebra automorphism
$\OmegabhmPi:\defUbhmPi\to\defUbhmPi$ with
\begin{equation*}
\begin{array}{l}
\OmegabhmPi(\defacK_\lambda\defacL_\mu)=\defacK_{-\lambda}\defacL_{-\mu}
\quad(\lambda,\mu\in\defVZ), \\
\OmegabhmPi(\defacE_i)=\defacF_i\defacL_{-\defPi(i)},\,\,
\OmegabhmPi(\defacF_i)=\defacK_{-\defPi(i)}\defacE_i
\quad(i\in\defI).
\end{array}
\end{equation*}

Define the bicharacter $\bhmop:\defVZ\times\defVZ\to\bKt$ by
$\bhmop(x,y):=\bhm(y,x)$
($x$, $y\in\defVZ$).
By Lemma~\ref{lemma:univ},
we have the $\bK$-algebra isomorphism
$\defUpbhmopPi: \defU(\bhmop,\defPi)\to\defUbhmPi$ with
\begin{align}\label{e:chevalley}
\defUpbhmopPi(\defacK_\lambda\defacL_\mu)&=\defacK_\mu\defacL_\lambda,\,(\lambda,\mu\in\defVZ),\\
\defUpbhmopPi(\defacE_i)&=\defacF_i,\,
\defUpbhmopPi(\defacF_i)=\defacE_i
\,(i\in\defI),
\end{align}
which will be referred as the {\it Chevalley involution}. We also have the $\bK$-algebra anti-automorphism
$\GammabhmopPi:\defUbhmopPi\to\defUbhmPi$ with
\begin{equation*}
\GammabhmopPi(\defacK_\lambda\defacL_\mu)=\defacK_\mu\defacL_\lambda,\,(\lambda,\mu\in\defVZ),\,\,
\GammabhmopPi(\defacE_i)=\defacE_i,\,\GammabhmopPi(\defacF_i)=\defacF_i\,(i\in\defI).
\end{equation*}

For $\defU=\defUbhmPi$, we have
the following equations. See \cite{Hec10} for the
notation ${{m}\choose{r}}_{q_{ii}}$ and ${{k}\choose{r}}_{q_{ii}}$.
\begin{equation*}
\begin{array}{lcl}
\defacE_i^k\defacF_i^m &
=& \sum_{r=0}^{{\mathrm{min}}\{k,m \}}(r)_{\tq_{ii}}!{m\choose r}_{\tq_{ii}}
{k\choose r}_{\tq_{ii}}
\tq_{ii}^{{\frac {r(-2k+r+1)} 2}} \\
& & \quad\quad\quad\quad\quad\cdot(\prod_{s=0}^{r-1}(-\defacK_{\defPi(i)}+\tq_{ii}^{-m+k+s}\defacL_{\defPi(i)}))
\defacF_i^{m-r}\defacE_i^{k-r},
\end{array}
\end{equation*} whence
\begin{equation*}
\begin{array}{l}
\mbox{$[\defacE_i,\defacF_i^m]$}=(m)_{\tq_{ii}}(-\defacK_{\defPi(i)}+\tq_{ii}^{-m+1}\defacL_{\defPi(i)})\defacF_i^{m-1}, \\
\mbox{$[\defacF_i,\defacE_i^m]$}=(m)_{\tq_{ii}}(-\defacL_{\defPi(i)}+\tq_{ii}^{-m+1}\defacK_{\defPi(i)})\defacE_i^{m-1}.
\end{array}
\end{equation*}
For $m\in\bZp$ and $i$, $j\in\defI$ with $i\ne j$,
define $\defacE_{m,i,j}$, $\defacEc_{m,i,j}\in\defUp_{\defPi(j)+m\defPi(i)}$, $\defacF_{m,i,j}$
$\defacF_{m,i,j}\in\defUm_{-\defPi(j)-m\defPi(i)}$
inductively by
$\defacE_{0,i,j}:=\defacEc_{0,i,j}:=\defacE_j$,
$\defacF_{0,i,j}:=\defacEc_{0,i,j}:=\defacF_j$, and
\begin{equation}\label{eqn:EFp}
\begin{array}{l}
\defacE_{m+1,i,j}:=\defacE_i\defacE_{m,i,j}
-\tq_{ii}^m\tq_{ij}\defacE_{m,i,j}\defacE_i,\,
\defacEc_{m+1,i,j}:=\defacEc_{m,i,j}\defacE_i
-\tq_{ii}^m\tq_{ji}\defacE_i\defacEc_{m,i,j},
\\
\defacF_{m+1,i,j}:=\defacF_i\defacF_{m,i,j}
-\tq_{ii}^m\tq_{ji}\defacF_{m,i,j}\defacF_i,\,
\defacFc_{m+1,i,j}:=\defacFc_{m,i,j}\defacF_i
-\tq_{ii}^m\tq_{ij}\defacF_i\defacFc_{m,i,j}.
\end{array}
\end{equation}
We have
\begin{equation}\label{eqn:defUpEmij}
\begin{array}{l}
\defUpbhmopPi(\defacE_{m,i,j})=\defacF_{m,i,j}, \quad
\defUpbhmopPi(\defacF_{m,i,j})=\defacE_{m,i,j}, \\
\defUpbhmopPi(\defacEc_{m,i,j})=\defacFc_{m,i,j}, \quad
\defUpbhmopPi(\defacFc_{m,i,j})=\defacEc_{m,i,j}.
\end{array}
\end{equation}
We have
\begin{equation}\label{eqn:tUbiddtwo}
\begin{array}{l}
\mbox{$[\defacE_i,\defacF_{m,i,j}]$}=-(m)_{\tq_{ii}}
(m;\tq_{ii},\tq_{ij}\tq_{ji})\defacK_{\defPi(i)}\defacF_{m-1,i,j},
\\
\mbox{$[\defacF_i,\defacE_{m,i,j}]$}=-(m)_{\tq_{ii}}
(m;\tq_{ii},\tq_{ij}\tq_{ji})\defacL_{\defPi(i)}\defacE_{m-1,i,j},
\\
\mbox{$[\defacE_j, \defacF_{m,i,j}]$}=
(m;\tq_{ii},\tq_{ij}\tq_{ji})!\defacF_i^m\defacL_{\defPi(j)},
\\
\mbox{$[\defacF_j, \defacE_{m,i,j}]$}=
(m;\tq_{ii},\tq_{ij}\tq_{ji})!\defacE_i^m\defacK_{\defPi(j)},
\\
\mbox{$[\defacE_{m,i,j},\defacF_{m,i,j}]$}=(m)_{\tq_{ii}}!
(m;\tq_{ii},\tq_{ij}\tq_{ji})!(-\defacK_{\defPi(j)+m\defPi(i)}+\defacL_{\defPi(j)+m\defPi(i)}),
\\
\mbox{$[\defacE_{m,i,j},\defacF_{m^\prime,i,k}]$}=0
\quad (m^\prime\in\bZp,k\in\defI\setminus\{i,j\}).
\end{array}
\end{equation} We have
\begin{equation}\label{eqn:OPiEFrij}
\begin{array}{l}
\OmegabhmPi(\defacE_{r,i,j})=\tq_{ii}^{-{\frac {r(r-1)} {2}}}\tq_{ji}^{-r}
\defacF_{r,i,j}\defacL_{-\defPi(j)-r\defPi(i)}, \\
\OmegabhmPi(\defacEc_{r,i,j})=\tq_{ii}^{-{\frac {r(r-1)} {2}}}\tq_{ij}^{-r}
\defacFc_{r,i,j}\defacL_{-\defPi(j)-r\defPi(i)}, \\
\OmegabhmPi(\defacF_{r,i,j})=\tq_{ii}^{\frac {r(r-1)} {2}}\tq_{ji}^r
\defacK_{-\defPi(j)-r\defPi(i)}\defacE_{r,i,j}, \\
\OmegabhmPi(\defacFc_{r,i,j})=\tq_{ii}^{\frac {r(r-1)} {2}}\tq_{ij}^r
\defacK_{-\defPi(j)-r\defPi(i)}\defacEc_{r,i,j}.
\end{array}
\end{equation} We have
\begin{equation}\label{eqn:XiPiEFrij}
\begin{array}{l}
\GammabhmopPi\GammabhmPi=\rmid_\defUbhmPi,\,\GammabhmPi\GammabhmopPi=\rmid_\defUbhmopPi,  \\
\GammabhmopPi(\defacE_{r,i,j})=\defacEc_{r,i,j},\,
\GammabhmopPi(\defacF_{r,i,j})=\defacFc_{r,i,j}.
\end{array}
\end{equation}

\subsection{Kharchenko's Poincar{\'{e}}-Birkohoff-Witt theorem and Heckenberger's Lusztig~isomorphisms}
We recall the following theorem by Kharchenko. We also introduce some notations.

\begin{theorem}\label{theorem:KharchenkoPBW}
{\rm{(\cite{Kha99})}}
Let $\bhm$ be a bicharacter and $\defU=\defUbhmPi$ be the generalized quantum group.
Then there exists a unique pair $(\defRbhmpip,\defvphibhmpip)$ of a
subset $\defRbhmpip$ of $\defVZPip$ and a map $\defvphibhmpip:\defRbhmpip\to\bN$
satisfying the following:

Let $X:=\{(\al,t)\in\defRbhmpip\times\bN|t\in\bJ_{1,\defvphibhmpip(\al)}\}$.
Define the map $z:X\to\defRbhmpip$ by $z(\al,t):=\al$.
Let $Y$ be the set of maps $y:X\to\bZp$ such that
$|\{x\in X|y(x)\geq 1\}|<\infty$ and
$(y(x))_{\bhm(z(x),z(x))}!\ne 0$ for all
$x\in X$. Then
\newline\newline
${\rm{(*)}}$\quad $\forall\lambda\in\defVZPip$,
$\dim\defUp_\lambda=|\{y\in Y|\sum_{x\in X}y(x)z(x)=\lambda\}|$.
\newline\newline Moreover, letting $\tq_{i^\prime j^\prime}:=\bhm(\defPi(i^\prime),\defPi(j^\prime))$
$(i^\prime, j^\prime\in\defI)$, for $i$, $j\in\defI$ with $i\ne j$, we have
\begin{equation}\label{eqn:goodbhmcond}
\{t\in\bZp|\defPi(j)+t\defPi(i)\in\defRbhmpip\}=\{t\in\bZp|(t)_{\tq_{ii}}!(t;\tq_{ii},\tq_{ij}\tq_{ji})!\ne 0\}.
\end{equation}
\end{theorem}

Let $\defRbhmpi:=\defRbhmpip\cup(-\defRbhmpip)$.
Define the map $\defvphibhmpi:\defRbhmpi\to\bN$
by $\defvphibhmpi(-\al):=\defvphibhmpi(\al):=\defvphibhmpip(\al)$
($\al\in\defRbhmpip$).
We say that $\bhm$ is {\it{$(\defPi,i)$-good}} if it satisfies the condition that
for every $j\in\defI\setminus\{i\}$,
there exists $\Nbhmpi_{ij}\in\bZp$ such that
$\defPi(j)+\Nbhmpi_{ij}\defPi(i)\in\defRbhmpip$ and
$\defPi(j)+(\Nbhmpi_{ij}+1)\defPi(i)\notin\defRbhmpip$.
Let $\Nbhmpi_{ii}:=-2$.

Let $i\in\defI$ and $\bhm$ be a $(\defPi,i)$-good bicharacter.
Define the map $\tauidefPi:\defI\to\defVZ$
by $\tauidefPi(j):=\defPi(j)+\Nbhmpi_{ij}\defPi(i)$
($j\in\defI$). By \eqref{eqn:goodbhmcond},
we see that $\tauidefPi$ is a $(\defPi,i)$-good bicharacter, and
\begin{equation}\label{eqn:Nbhminv}
\NbhmdtiPi_{ij}=\Nbhmpi_{ij}\quad(j\in\defI),
\end{equation} whence $(\deftaubhm_i)^2\defPi=\defPi$.

\begin{theorem}\label{theorem:HecIso}
{\rm{(\cite{Hec10})}}
Let $a:\defI\to\bKt$ be a map.
Let $i\in\defI$. Let $\bhm$ be a $(\defPi,i)$-good bicharacter.
Then there exists a
$\bK$-algebra isomorphism
$\LusT^{\bhm,\tauidefPi,a}_i:\defUbhmtauiPi\to\defUbhmPi$ such that
\begin{equation}\label{eqn:dfLusT}
\begin{array}{l}
\LusT_i(\defacK_\lambda)=\defacK_\lambda,\quad \LusT_i(\defacL_\lambda)=\defacL_\lambda
\quad (\lambda\in\defVZ),
\\
\LusT_i(\defacE_i)=a(i)\defacF_i\defacL_{-\defPi(i)},\quad
\LusT_i(\defacF_i)={\frac 1 {a(i)}}\defacK_{-\defPi(i)}\defacE_i,\\
\LusT_i(\defacE_j)=a(j)\defacE_{\Nbhmpi_{ij},i,j} \quad (j\in\defI\setminus\{i\})\\
\LusT_i(\defacF_j)={\frac 1
{a(j)(\Nbhmpi_{ij})_{\tq_{ii}}!(\Nbhmpi_{ij};\tq_{ii},\tq_{ij}\tq_{ji})!}}
\defacF_{\Nbhmpi_{ij},i,j} \quad (j\in\defI\setminus\{i\}),
\end{array}
\end{equation}
where $\LusT_i:=\LusT^{\bhm,\tauidefPi,a}_i$ and $\tq_{i^\prime j^\prime}:=\bhm(\defPi(i^\prime),\defPi(j^\prime))$
$(i^\prime,\,j^\prime\in\defI)$. Moreover we have
\begin{equation}\label{eqn:HecIsoEq}
\defRtaubhmipi=\defRbhmpi,\,\,\defRtaubhmipip=(\defRbhmpip\setminus\{\defPi(i)\})\cup\{-\defPi(i)\},\,\,
\defvphibhmpi=\defvphibhmtauipi.
\end{equation}
\end{theorem}

Let $i$, $j\in\defI$ be such that $i\ne j$.
Let $\bhm$ be a $(\defPi,i)$-good and $(\defPi,j)$-good bicharacters.
We say that $\bhm$ is {\it{a $(\defPi,i,j)$-good bicharacter}}
if $\deftaubhm_{i_1}\cdots\deftaubhm_{i_m}\defPi$
can be defined for all $m\in\bN$ and
all $i_t\in\{i,j\}$
($t\in\bJ_{1,m}$).

\begin{equation}\label{eqn:defmij}
\mbox{Let $\mbhmpi_{ij}:=|\defRbhmpip\cap(\bZ\defPi(i)\oplus\bZ\defPi(j))|(\in\bJ_{2,\infty}\cup\{\infty\})$.}
\end{equation}

By \eqref{eqn:HecIsoEq}, using the same argument as that for \cite[Lemma~4]{HY08}, we have the following result.
\begin{lemma}\label{lemma:HecIRktetwo}
Let $\bhm$ be a $(\defPi,i,j)$-good bicharacter.
Let $X:=\bZ\defPi(i)\oplus\bZ\defPi(j)$.
Let $m:=\mbhmpi_{ij}$.
Let $i_{2y-1}:=i$, $i_{2y}:=j$ $(y\in\bN)$.
Let $\defPi_t:=\deftaubhm_{i_t}\deftaubhm_{i_{t-1}}\cdots\deftaubhm_{i_1}\defPi$
$(t\in\bN)$. Then for
$k\in\bJ_{1,m}$, we have
\begin{equation*}
|\defRbhmpip\cap(-R_\bhm^{\defPi_k,+})|=|\defRbhmpip\cap(-R_\bhm^{\defPi_k,+})\cap X|=k.
\end{equation*}
Moreover if $m<\infty$, we have
\begin{equation*}
\defPi_m=\deftaubhm_{i_{m+1}}\deftaubhm_{i_m}\cdots\deftaubhm_{i_2}\defPi,
\end{equation*} and
\begin{equation*}
\defRbhmpip\cap X=\{\defPi(i_1)=i,\defPi_1(i_2),\ldots,\defPi_{m-2}(i_{m-1}),
\defPi_{m-1}(i_m)=j\}.
\end{equation*}
Furthermore $|\{\defPi_t(i_{t+1})|t\in\bN\}|=\infty$ if $m=\infty$.
\end{lemma}

Let $\bhm$ be a $(\defPi,i,j)$-good bicharacter.
We say that $\bhm$ is {\it{a $(\defPi,i,j)$-good finite bicharacter}}
if $m$ in Lemma~\ref{lemma:HecIRktetwo} is finite.

\begin{theorem}\label{theorem:HecIfund} {\rm{(\cite{Hec10})}}
Let $i, j \in I, i\neq j$ and $\bhm$ be a $(\defPi,i,j)$-good finite bicharacter with $m:=\mbhmpi_{ij} <\infty$.
For $t\in\bN$, let $\defPi^\prime_t:=\deftaubhm_{i_{t+1}}\deftaubhm_{i_t}\cdots\deftaubhm_{i_2}\defPi$
{\rm{(}}$\mbox{$i_k = i$ {\rm{(}}{\rm resp.} $j${\rm{)}}}$ if $k$ is odd {\rm{(}}resp. even{\rm{)}}{\rm{)}},
 $a_t$, $a^\prime_t:\defI\to\bKt$ be maps,
${\dot{\LusT}}_{(t)}:=\LusT^{\bhm,\defPi_1,a_1}_{i_1}\cdots\LusT^{\bhm,\defPi_t,a_t}_{i_t}$
and ${\ddot{\LusT}}_{(t)}:=\LusT^{\bhm,\defPi^\prime_1,a^\prime_1}_{i_2}\cdots\LusT^{\bhm,\defPi^\prime_t,a^\prime_t}_{i_{t+1}}$.
Then there exists a map $b:\defI\to\bKt$ such that
\begin{equation*}
{\dot{\LusT}}_{(m)}(\defacE_k)=b(k){\ddot{\LusT}}_{(m)}(\defacE_k),\,\,
{\dot{\LusT}}_{(m)}(\defacF_k)={\frac 1 {b(k)}}{\ddot{\LusT}}_{(m)}(\defacF_k)
\quad (k\in\defI).
\end{equation*}
Moreover there exists $z\in\bKt$ such that
\begin{equation*}
{\dot{\LusT}}_{(m-1)}(\defacE_{i_m})=z\defacE_j,\,\,{\dot{\LusT}}_{(m-1)}(\defacF_{i_m})=
z^{-1}\defacF_j.
\end{equation*}
\end{theorem}

\subsection{Strict Heckenberger's Lusztig isomorphisms}

\begin{lemma}\label{lemma:plusminusone}
Let $i, j \in I, i\neq j$ and $\bhm$ be a $(\defPi,i,j)$-good finite bicharacter with $m:=\mbhmpi_{ij} <\infty$. Also
$b:\defI\to\bKt$ be the map and $z\in\bKt$  in Theorem~{\rm{\ref{theorem:HecIfund}}}. Then the following statements hold.
 \newline\newline
{\rm{(1)}} If $\OmegabhmPi{\dot{\LusT}}_{(m)}={\dot{\LusT}}_{(m)}\defOmega^{\bhm,\defPi_m}$,
then $b(k)\in\{1,-1\}$ for all $k\in\defI$.\newline
{\rm{(2)}} If $\OmegabhmPi{\dot{\LusT}}_{(m-1)}={\dot{\LusT}}_{(m-1)}\defOmega^{\bhm,\defPi_{m-1}}$,
then $z\in\{1,-1\}$.
\end{lemma}
{\it{Proof.}} (1) Let $k\in\defI$. We have
\begin{equation}
\begin{array}{l}
\defacF_k=b(k){\ddot{\LusT}}_{(m)}^{-1}{\dot{\LusT}}_{(m)}(\defacF_k)
=b(k){\ddot{\LusT}}_{(m)}^{-1}{\dot{\LusT}}_{(m)}\defOmega^{\bhm,\defPi_m}(\defacE_k
\defacL_{\defPi(k)}) \\
\quad =b(k)\OmegabhmPi{\ddot{\LusT}}_{(m)}^{-1}{\dot{\LusT}}_{(m)}(\defacE_k
\defacL_{\defPi(k)})
=b(k)^2\OmegabhmPi(\defacE_k\defacL_{\defPi(k)})
=b(k)^2\defacF_k,
\end{array}
\end{equation} whence $b(k)^2=1$, so $b(k)\in\{1,-1\}$. \newline\par
(2) This can be proved similarly to (1).
\hfill $\Box$.

\begin{lemma}\label{lemma:OmegaT}
Let $a:\defI\to\bKt$ be a map.
Let $i\in\defI$ and $\bhm$ be a $(\defPi,i)$-good bicharacter.
Then
\begin{equation*}
\OmegabhmPi\LusT^{\bhm,\tauidefPi,a}_i
=\LusT^{\bhm,\tauidefPi,a}_i\defOmega^{\bhm,\tauidefPi}
\end{equation*} if and only if
\begin{equation}\label{eqn:dfOmegaT}
a(i)^2=1\quad\mbox{and}\quad a(j)^2={\frac {\tq_{ii}^{\frac {\Nbhmpi_{ij}(\Nbhmpi_{ij}-1)} 2}\tq_{ji}^{\Nbhmpi_{ij}}}
{(\Nbhmpi_{ij})_{\tq_{ii}}!(\Nbhmpi_{ij};\tq_{ii},\tq_{ij}\tq_{ji})!}}\,\,
(j\in\defI\setminus\{i\}).
\end{equation}
\end{lemma}
{\it{Proof.}} The claim follows from Theorem~\ref{theorem:HecIso} and \eqref{eqn:OPiEFrij}
\hfill $\Box$

\begin{lemma}\label{lemma:sltwoac}
Let $G$ be the $\bK$-algebra {\rm{(}}with $1${\rm{)}} defined with
the generators $X$, $Y$, $Z$ and the relations
$[Z,X]=2X$, $[Z,Y]=-2Y$, $[X,Y]=Z$.
{\rm{(}}Namely  $G$ is isomorphic to the universal enveloping algebra
of ${\mathrm{sl}}_2(\bK)$.{\rm{)}}
Let $k\in\bZp$. Let $\Gamma$ be the
$(k+1)$-dimensional left $G$-module
with the $\bK$-basis $\{\gamma_r|r\in\bJ_{0,k}\}$
such that $Z\gamma_r:=(k-2x)\gamma_r$,
$Y\gamma_r:=\gamma_{r+1}$, $X\gamma_r:=x(k-r+1)\gamma_{r-1}$,
where $\gamma_{-1}:=\gamma_{k+1}:=0$.
Let $a\in\bKt$. Define $\smallT\in{\rm{End}}_\bK(\Gamma)$
by $\smallT(v):=\exp(aX)\exp(-a^{-1}Y)\exp(aX)\cdot v$
$(v\in\Gamma)$. Then
\begin{equation*}\label{eqn:gYrmo}
\smallT(Y^r\gamma_0)={\frac{(-1)^{k-r}a^{2r-k}r!} {(k-r)!}}Y^{k-r}\gamma_0, \quad
\smallT(X^r\gamma_k)={\frac{(-1)^ra^{k-2r}r!} {(k-r)!}}X^{k-r}\gamma_k.
\end{equation*}
\end{lemma}
{\it{Proof.}} Let ${\bar \Gamma}$, ${\bar \gamma}_y$ ($y\in\bJ_{0,1}$),
${\bar \smallT}$ be $\Gamma$, $\gamma_y$, $\smallT$ respectively for $k=1$.
Then ${\bar \smallT}({\bar \gamma}_0)=(-a^{-1}){\bar \gamma}_1$ and ${\bar \smallT}({\bar \gamma}_1)=a{\bar \gamma}_0$.
Regard ${\bar \Gamma}^{\otimes k}$ as the ($k$-fold tensor) $G$-module
in a standard way. Let $g:M\to{\bar M}^{\otimes k}$ be the $G$-module monomorphism $g:M\to{\bar \Gamma}^{\otimes k}$
such that $g(\gamma_0)={\bar \gamma}_0^{\otimes k}$.
Then
\begin{equation*}
g(\gamma_x)=g(Y^x\gamma_0)=Y^x\cdot g(\gamma_0)=r!\sum_{|\{x\in\bJ_{1,k}|i_x=1\}|=x}{\bar \gamma}_{i_1}\otimes\cdots\otimes{\bar \gamma}_{i_1}.
\end{equation*}
Note $g(\smallT(v))={\bar \smallT}^{\otimes k}(g(v))$ ($v\in\Gamma$).
Then
\begin{equation*}
\begin{array}{l}
g(\smallT(\gamma_r))=g(t(Y^r\gamma_0))={\bar \smallT}^{\otimes k}(g(Y^r\gamma_0))\\
\quad ={\frac{(-a^{-1})^{k-r}a^rr!} {(k-r)!}}g(Y^{k-r}\gamma_0)
=g({\frac{(-1)^{k-r}a^{2r-k}r!} {(k-r)!}}\gamma_{k-r}).
\end{array}
\end{equation*}
Thus we can see the claim.
\hfill $\Box$

\begin{lemma}\label{lemma:ranktwot}
Let $a_{11}:=a_{22}:=2(\in\bK)$. Let
$a_{12}$, $a_{21}\in\bK$ be such that
$(a_{12},a_{21})\in\{(0,0),(-1,-1),(-2,-1),(-3,-1)\}$.
Let $m\in\bZp$ be $2$ {\rm{(}}{\rm{resp.}} $3$, {\rm{resp.}} $4$, {\rm{resp.}} $6${\rm{)}}
if $a_{12}$ is $0$ {\rm{(}}{\rm{resp.}} $-1$, {\rm{resp.}} $-2$, {\rm{resp.}} $-3${\rm{)}}.
Let $i_{2x-1}:=1$, $i_{2x}:=2$ $(x\in\bN)$.
Let $G$ be a $\bK$-algebra {\rm{(}}with $1${\rm{)}}
satisfying the following conditions {\rm{(i)}} and {\rm{(ii)}}.
\newline\newline
{\rm{(i)}}
There exist $X_i$, $Y_i$, $Z_i\in G$
$(i\in\bJ_{1,2})$ such that
\begin{equation*}
\begin{array}{l}
[Z_1,Z_2]=0,\,[Z_i,X_j]=a_{ij}X_j,\,[Z_i,Y_j]=-a_{ij}Y_j,\,[X_i,Y_j]=\delta_{ij}Z_i, \\
(\rmad X_{i^\prime})^{1-a_{i^\prime j^\prime}}(X_{j^\prime})
=(\rmad Y_{i^\prime})^{1-a_{i^\prime j^\prime}}(Y_{j^\prime})=0\,\,(i^\prime\ne j^\prime).
\end{array}
\end{equation*}
\newline
{\rm{(ii)}} For $\Re\in G$ and $\Im\in\{X_i,Y_i|i\in\bJ_{1,2}\}$, there exists $r\in\bN$
such that $(\rmad \Im)^r(\Re)=0$.
\newline\newline
For $i\in\bJ_{1,2}$, let $b_i\in\bKt$ $(i\in\bJ_{1,2})$, and
define the $\bK$-algebra automorphism $\smallT_i$ by
$\smallT_i(\Re):=\exp(b_i\rmad X_i)\exp(-b_i^{-1}\rmad Y_i)\exp(b_i\rmad X_i)(\Re)$
$(\Re\in G)$.
Then we have
\begin{equation}\label{eqn:rktwota}
\begin{array}{l}
\smallT_i(X_i)=-b_i^2Y_i,\,\,\smallT_i(Y_i)=-b_i^{-2}X_i, \\
\smallT_i(X_j)={\frac {b_i^{-a_{ij}}} {(-a_{ij})!}}(\rmad X_i)^{-a_{ij}}(X_j),\,\,
\smallT_i(Y_j)={\frac {(-b_i)^{a_{ij}}} {(-a_{ij})!}}(\rmad Y_i)^{-a_{ij}}(Y_j) \quad (i\ne j),
\end{array}
\end{equation}
\begin{equation}\label{eqn:rktwotb}
\begin{array}{l}
\smallT_{i_t}\smallT_{i_{t+1}}\cdots\smallT_{i_{t+m-2}}(X_{i_{t+m-1}})=X_{i_{t-1}}, \\
\smallT_{i_t}\smallT_{i_{t+1}}\cdots\smallT_{i_{t+m-2}}(Y_{i_{t+m-1}})=Y_{i_{t-1}} \\
(t\in\bZ),
\end{array}
\end{equation} and
\begin{equation}\label{eqn:rktwotc}
\smallT_{i_1}\smallT_{i_2}\cdots\smallT_{i_m}
=\smallT_{i_2}\smallT_{i_3}\cdots\smallT_{i_{m+1}}.
\end{equation}
\end{lemma}
{\it{Proof.}} We can see \eqref{eqn:rktwota} by Lemma~\ref{lemma:sltwoac}.
As for \eqref{eqn:rktwotb},
for example, if $a_{12}=-3$, by Lemma~\ref{lemma:sltwoac}, we have
\begin{equation*}
\begin{array}{l}
\smallT_1\smallT_2\smallT_1\smallT_2\smallT_1(X_2) \\
\quad = \smallT_1\smallT_2\smallT_1\smallT_2({\frac {b_1^3} {3!}}[X_1,[X_1,[X_1,X_2]]]) \\
\quad = -{\frac {b_1^3} {3!}}\smallT_1\smallT_2\smallT_1\smallT_2([X_1,[X_1,[X_2,X_1]]]) \\
\quad = -{\frac {b_1^3} {3!}}\smallT_1\smallT_2\smallT_1(
[b_2[X_2,X_1],[b_2[X_2,X_1],(-1)b_2^{-1}X_1]]) \\
\quad = -{\frac {b_1^3b_2} {3!}}\smallT_1\smallT_2(
[(-1){\frac {b_1} {2!}}[X_1,[X_1,X_2]],{\frac {2!} {b_1}}[X_1,X_2]]) \\
\quad ={\frac {b_1^3b_2} {3!}}\smallT_1\smallT_2(
[[X_1,[X_2,X_1]],[X_2,X_1]]) \\
\quad ={\frac {b_1^3b_2} {3!}}\smallT_1(
[[b_2[X_2,X_1],(-1)b_2^{-1}X_1],(-1)b_2^{-1}X_1]) \\
\quad =-{\frac {b_1^3} {3!}}\smallT_1([X_1,[X_1,[X_1,X_2]]]) \\
\quad =(-1){\frac {b_1^3} {3!}}\cdot (-1){\frac {3!} {b_1^3}}X_2 \\
\quad =X_2.
\end{array}
\end{equation*}
Now we show \eqref{eqn:rktwotc}. Let $\smallT^\prime:=\smallT_{i_1}\smallT_{i_2}\cdots\smallT_{i_{m-1}}$.
By \eqref{eqn:rktwotb}, for $\Re\in G$, we see
\begin{equation*}
\begin{array}{l}
\smallT^\prime\smallT_{i_m}(\smallT^\prime)^{-1}(\Re) \\
\quad =\exp(b_i\rmad \smallT^\prime(X_{i_m}))\exp(-b_i^{-1}\rmad \smallT^\prime(Y_{i_m}))
\exp(b_i\rmad \smallT^\prime(X_{i_m}))(\Re) \\
\quad =\smallT_{i_2}(\Re).
\end{array}
\end{equation*}
This completes the proof. \hfill $\Box$
\newline\par
\begin{definition}\label{definition:admbhm}
Let $\bhm:\defVZ\times\defVZ\to\bKt$ be a bi-character, $\defPi:\defI\to\defVZ$ be a map such that
$\defPi(\defI)$ is a $\bZ$-base of $\defVZ$ and $\tq_{ij}:=\bhm(\defPi(i),\defPi(j))$ for all $i,j \in I$.
Let $A=[a_{ij}]_{ij\in\defI}$ be a symmetrizable generalized Cartan matrix.
Let $d_i\in\bN$ be such that $d_ia_{ij}=d_ja_{ji}$.
($i$, $j\in\defI$).
Let $\dotq\in\bKt$ be such that
$\dotq^r\ne 1$ for all $r\in\bN$.
Let $\sqrtbhm:\defVZ\times\defVZ\to\bKt$ be a bicharacter and $\dotq_{ij}:=\sqrtbhm(\defPi(i),\defPi(j))$ ($i$, $j\in\defI$).
Assume that $\dotq_{ii}=\dotq^{2d_i}$ $(i\in\defI)$,
$\dotq_{ij}\dotq_{ji}=\dotq^{2d_ia_{ij}}$ $(i,j\in\defI)$,
$q_{ij}:=\dotq_{ij}^2$ for all $i$, $j\in\defI$.
Also assume that for every $k\in\defI$, there exists $\defTheta(\tq_{kk}-1)\in\bKt$ such that
$\defTheta(\tq_{kk}-1)^2=\tq_{kk}-1$.
Then we say that such $\bhm$ is {\it{a $(\defPi,A)$-admissible bicharacter}}.
If $A$ is the Cartan matrix of a finite-dimensional complex Lie algebra
(i.e., $A$ is a symmetrizable generalized Cartan matrix of finite-type),
we call $\defUbhmPi$ {\it{a finite-type multiparameter quantum group.}}
\end{definition}

For a $(\defPi,A)$-admissible bicharacter $\bhm$, $\defUbhmPi$ is presented by the generators given by
Definition~\ref{definition:defofGQG}~(i)
and the relations composed of those of Definition~\ref{definition:defofGQG}~(ii)
and $\defacE_{1-a_{ij},i,j}=\defacF_{1-a_{ij},i,j}=0$ ($i,j\in\defI$, $i\ne j$),
which is well-known and can be proved by a standard argument along with Theorem~\ref{theorem:LusPBW} below.

\begin{lemma} \label{lemma:lmfundaTh}
Let $\bhm$ be a $(\defPi,A=[a_{ij}]_{ij\in\defI})$-admissible bicharacter.
\newline\newline {\rm{(1)}}  Then $\bhm$ is $(\defPi,i)$-good bi-character for every $i\in\defI$
and $\Nbhmpi_{ij}=-a_{ij}$ for all $i$, $j\in\defI$.
\newline {\rm{(2)}}
For $i\in\defI$ let
$\ddotq_{jk}:=\sqrtbhm(\tauidefPi(j),\tauidefPi(k))$
$(j,k\in\defI)$.
Then we have
\begin{equation*}
\ddotq_{jj}=\dotq_{jj}
\,\,(j\in\defI),\,\,
\ddotq_{jk}\ddotq_{kj}
=\dotq_{jk}\dotq_{kj}\,\,(j,k\in\defI).
\end{equation*}
In particular, $\bhm$ is $(\tauidefPi,A)$-admissible.
\newline {\rm{(3)}} For $i$, $j\in\defI$ with $i\ne j$,
$\bhm$ is $(\defPi,i,j)$-finite if and only if
$a_{ij}a_{ji}\in\bJ_{0,3}$. Moreover
\begin{equation}\label{eqn:prmij}
\mbhmpi_{ij}=
\left\{\begin{array}{ll}
2 & \mbox{if $a_{ij}a_{ji}=0$,} \\
3 & \mbox{if $a_{ij}a_{ji}=1$,} \\
4 & \mbox{if $a_{ij}a_{ji}=2$,} \\
6 & \mbox{if $a_{ij}a_{ji}=3$.}
\end{array}\right.
\end{equation}
\end{lemma}
{\it{Proof.}} Claim~(1) follows from \eqref{eqn:goodbhmcond}.
Claim~(2) can be proved directly.
Claim~(3) follows from Claims~(1) and (2) and Lemma~\ref{lemma:HecIRktetwo}. \hfill $\Box$
\newline\par

Let $\bhm$ be a $(\defPi,A=[a_{ij}]_{ij\in\defI})$-admissible bicharacter.
Let
$\dotq_{jk}$ be as in Definition~\ref{definition:admbhm}. Define the map
$\spa:\defI\to\bKt$ by
\begin{equation*}
\spa(j)=
\left\{\begin{array}{ll}
1 & \quad (j=i), \\
{\frac {\dotq_{ij}^{a_{ij}}}
{(-a_{ij})_{\tq_{ii}}!\defTheta(\tq_{ii}-1)^{-a_{ij}}}} & \quad(j\in\defI\setminus\{i\}).
\end{array}\right.
\end{equation*}
Let
\begin{equation*}
\LusTi:=\LusT^{\bhm,\tauidefPi,\spa}_i.
\end{equation*}
We see directly that $\spa$ satisfies \eqref{eqn:dfOmegaT}.
As for \eqref{eqn:dfLusT}, letting $\LusT_i:=\LusTi$, we have
\begin{equation*}
\begin{array}{l}
\LusT_i(\defacK_\lambda\defacL_\mu)=\defacK_\lambda\defacL_\mu \,
(\lambda,\mu\in\defVZ),\,\,
\LusT_i(\defacE_i)=\defacF_i\defacL_{-\defPi(i)},\,\LusT_i(\defacF_i)=\defacK_{-\defPi(i)}\defacE_i, \\
\LusT_i(\defacEc_{r,i,j})={\frac {(r)_{\tq_{ii}}!\dotq_{ij}^{a_{ij}+2r}\defTheta(\tq_{ii}-1)^{a_{ij}+2r}}
{(-a_{ij}-r)_{\tq_{ii}}!}}\defacE_{-a_{ij}-r,i,j},\\
\LusT_i(\defacFc_{r,i,j})={\frac {(r)_{\tq_{ii}}!\dotq_{ii}^{(a_{ij}+2r)(a_{ij}-1)}\dotq_{ij}^{-a_{ij}-2r}\defTheta(\tq_{ii}-1)^{a_{ij}+2r}}
{(-a_{ij}-r)_{\tq_{ii}}!}}\defacF_{-a_{ij}-r,i,j}
\\ (j\in\defI\setminus\{i\},\,r\in\bJ_{0,-a_{ij}}).
\end{array}
\end{equation*} We also have
\begin{equation} \label{eqn:InvLusTi}
(\LusTidefPi)^{-1}=\GammabhmopPi\LusTiop\GammabhmtauiPi.
\end{equation}

For $i\in\defI$, define the $\bK$-algebra automorphism
$\zetabhmPi_i:\defUbhmPi\to\defUbhmPi$ by
\begin{equation*}
\begin{array}{l}
\zetabhmPi_i(\defacK_\lambda\defacL_\mu):=\defacK_\lambda\defacL_\mu\,\,
(\lambda, \mu\in\defVZ), \\
\zetabhmPi_i(\defacE_j):={\frac 1 {\dotq_{ij}\dotq_{ji}}}\defacE_j,\,
\zetabhmPi_i(\defacF_j):=\dotq_{ij}\dotq_{ji}\defacF_j\,\,
(j\in\defI).
\end{array}
\end{equation*} Then we have
\begin{equation*}\label{eqn:Lustzata}
\LusTi\defUpsilon^{\bhmop,\tauidefPi}=\zetabhmPi_i\defUpbhmopPi
\LusTiop\quad(i\in\defI).
\end{equation*}
For $i\in\defI$, let
\begin{equation*}
\BdefacE_i:={\frac {\defacE_i} {\defTheta(\tq_{ii}-1)}},\quad
\BdefacF_i:=-{\frac {\defacF_i} {\defTheta(\tq_{ii}-1)}},\quad
\BdefacH_i:={\frac {\defacK_{\defPi(i)}-\defacL_{\defPi(i)}} {\tq_{ii}-1}}.
\end{equation*}
Note that
\begin{equation*}
[\BdefacE_i,\BdefacF_i]=\delta_{ij}\BdefacH_i,
\end{equation*} and
\begin{equation*}
\BdefacH_i\BdefacE_j=\tq_{ij}\BdefacE_j\BdefacH_i+\tq_{ij}^{-1}
{\frac {\tq_{ij}\tq_{ji}-1} {\tq_{ii}-1}}\BdefacE_j\defacL_{\defPi(i)}.
\end{equation*}

Let $\bO$ be
the $\bQ$-subalgebra of $\bK$ generated by
$\dotq^{\pm 1}_{ij}$, ${\frac 1 {(-a_{ij})_{\tq_{ii}}!}}$ for all $i,j\in\defI$.
Let $\defUO$ (resp. $\defUoO$, resp. $\defUpO$, resp. $\defUmO$)
be the $\bO$-subalgebra (with $1$) of $\defU=\defUbhmPi$ (resp. $\defUo$, resp. $\defUp$, resp. $\defUm$)
generated by
$\defacK_{\defPi(i)}^{\pm 1}$,  $\defacL_{\defPi(i)}^{\pm 1}$, $\BdefacH_i$, $\BdefacE_i$, $\BdefacF_i$
(resp. $\defacK_{\defPi(i)}^{\pm 1}$,  $\defacL_{\defPi(i)}^{\pm 1}$, $\BdefacH_i$,
resp. $\BdefacE_i$, resp. $\BdefacF_i$) for all $i\in\defI$.
We can see
$\defUO=\defUmO\otimes_\bO\defUoO\otimes_\bO\defUpO$.
We also see that the elements
\begin{equation*}
\prod_{i\in\defI}\defacK_{\defPi(i)}^{x(i)}(\defacK_{\defPi(i)}\defacL_{\defPi(i)})^{y(i)}
\BdefacH_i^{z(i)}
\end{equation*}
($x(i)\in\{0,1\}$, $z(i)\in\bZp$, $y(i)\in\bZ$) form
$\bO$-basis of $\defUoO$.

\begin{theorem} \label{theorem:fundaTh}
Let $i, j \in I, i\neq j$ and $\bhm$ be a $(\defPi,i,j)$-good finite bicharacter with $m:=\mbhmpi_{ij} <\infty$.
Assume that $\bhm$ is $(\defPi,A)$-admissible. Let $b:\defI\to\bKt$ and $z\in\bKt$ be as in Theorem~{\rm{\ref{theorem:HecIfund}}}.
Assume that $\LusT^{\bhm,\defPi_t,a_t}_{i_t}=\LusT^{\bhm,\defPi_t}_{i_t}$
and $\LusT^{\bhm,\defPi^\prime_t,a^\prime_t}_{i_{t+1}}=\LusT^{\bhm,\defPi^\prime_t}_{i_{t+1}}$
$(t\in\bN)$. Then we have
\begin{equation}\label{eqn:maineq}
b(k)=1\,\,(k\in\defI)\quad{and}\quad{z=1}.
\end{equation}
\end{theorem}
{\it{Proof.}}
We divide the proof into steps.

{\it{Step}}~1. {\it{Assume that $\dotq$ is transcendental over $\bQ$
and that $\dotq_{j^\prime i^\prime}=1$
for all $i^\prime$, $j^\prime\in\defI$ with $i^\prime < j^\prime$.}}
Then $\bO$ is the
$\bQ$-subalgebra of $\bK$ generated by
$\dotq^{\pm 1}$, ${\frac 1 {(-a_{i^\prime j^\prime})_{\tq_{i^\prime i^\prime}}!}}$
for all $i^\prime$, $j^\prime\in\defI$ with $i^\prime \ne j^\prime$.
So $\bO$ is a principal integral domain.

Consider the $\bQ$-algebra $(\defUO)_1:=\defUO/(\dotq-1)\defUO$.
Let $f:\defUO\to (\defUO)_1$ be the canonical map.
Note that for $k\in\defI$, we have $f(\defacL_{\defPi(k)})=f(\defacK_{\defPi(k)})$,
and we see that $f(K_{\defPi(k)})$ is a central element of $(\defUO)_1$.
So we can consider the quotient $\bQ$-algebra $(\defUO)_1/(f(\defacK_{\defPi(k)})-1)$.
By Lemmas~\ref{lemma:plusminusone}, \ref{lemma:OmegaT} and \ref{lemma:ranktwot},
we have \eqref{eqn:maineq}.

{\it{Step}}~2. Let $\defU^\prime$ denote the $\defU$ of Step~1, and
let $\ddotq$ denote $\dotq$ for $\defU^\prime$.
Then $\ddotq$ is transcendental over $\bQ$.
Assume that that $\dotq_{j^\prime i^\prime}=1$
for all $i^\prime$, $j^\prime\in\defI$ with $i^\prime < j^\prime$.
We use a specialization argument with $\ddotq\to\dotq$;
this $\dotq$ is the one for $\defU$ of this step.
We can obtain \eqref{eqn:maineq}
from Step~1 by considering
the $\bO$-subalgebra of $\defU^\prime$
generating by $\defacK_{\defPi(t)}^{\pm 1}$,  $\defacL_{\defPi(t)}^{\pm 1}$,
$\defacE_t$, $\defacF_t$ ($t\in\defI$)
and using Lemma~\ref{lemma:univ}.

{\it{Step}}~3. {\it{General cases.}}
Repeat the same arguments as in Step~2.
\hfill $\Box$

\subsection{Drinfeld pairing}

Let $\defU=\defUbhmPi$ be the generalized quantum group.
We regard $\defU=\defUbhmPi$ as a Hopf algebra $(\defU,\tHopfD,\tHopfS,\tHopfe)$
with
$\tHopfD(\defacK_\lambda)=\defacK_\lambda\otimes\defacK_\lambda$,
$\tHopfD(\defacL_\lambda)=\defacL_\lambda\otimes\defacL_\lambda$,
$\tHopfD(\defacE_i)=\defacE_i\otimes  1+\defacK_{\defPi(i)}\otimes\defacE_i$,
$\tHopfD(\defacF_i)=\defacF_i\otimes
\defacL_{\defPi(i)}+1\otimes\defacF_i$,
$\tHopfS(\defacK_\lambda)=\defacK_{-\lambda}$,
$\tHopfS(\defacL_\lambda)=\defacL_{-\lambda}$,
$\tHopfS(\defacE_i)=-\defacK_{-\defPi(i)}\defacE_i$,
$\tHopfS(\defacF_i)=-\defacF_i\defacL_{-\defPi(i)}$,
$\tHopfe(\defacK_\lambda)=\tHopfe(\defacL_\lambda)=1$,
and $\tHopfe(\defacE_i)=\tHopfe(\defacE_i)=0$.

For $i$, $j\in\defI$ with $i\ne j$ and $r\in\bZp$, if $E_{
r,i,j}\ne 0$, we have
\begin{equation*}
\tHopfD(\defacE_{r,i,j})=\defacE_{r,i,j}\otimes 1
+\sum_{k=0}^r
{\frac {(r)_{\tq_{ii}}!(k;\tq_{ii},\tq_{ii}^{r-k}\tq_{ij}\tq_{ji})!} {(k)_{\tq_{ii}}!(r-k)_{\tq_{ii}}!}}
\defacE_i^k\defacK_{(r-k)\defPi(i)+\defPi(j)}\otimes\defacE_{r-k,i,j}.
\end{equation*}

Let $\tBp=\tBpbhmPi:=\oplus_{\lambda\in\defVZ}\defUp\defacK_\lambda$,
and $\tBm=\tBmbhmPi:=\oplus_{\lambda\in\defVZ}\defUm\defacL_\lambda$.
Then $\defU=\rmSpan_\bK(\tBm\tBp)=\rmSpan_\bK(\tBp\tBm)$ and in a
standard way (see \cite{Dr86}), we have a bilinear form
$\tvt=\tvtbhmPi:\tBp\times\tBm\to\bK$ with the following
properties:
\begin{equation}\label{eqn:bprtoftB}
\begin{array}{l}
\tvt(\defacK_\lambda,\defacL_\mu)=\bhm(\lambda,\mu),
\tvt(\defacE_i,\defacF_j)=\delta_{ij},
\tvt(\defacK_\lambda,\defacF_j)=\tvt(\defacE_i,\defacL_\lambda)=0, \\
\tvt(X^+Y^+,X^-)=
\sum_{k^-}\tvt(X^+,(X^-)^{(2)}_{k^-})\tvt(Y^+,(X^-)^{(1)}_{k^-}),\\
\tvt(X^+,X^-Y^-)=
\sum_{k^+}\tvt((X^+)^{(1)}_{k^+},X^-)\tvt((X^+)^{(2)}_{k^+},Y^-),\\
\tvt(\tHopfS(X^+),X^-)=\tvt(X^+,\tHopfS^{-1}(X^-)),\\
\tvt(X^+,1)=\tHopfe(X^+),
\tvt(1,X^-)=\tHopfe(X^-), \\
X^-X^+  =\sum_{r^+,r^-}
\tvt((X^+)^{\prime,(1)}_{r^+}, \tHopfS((X^-)^{\prime,(1)}_{r^-}))
\tvt((X^+)^{\prime,(3)}_{r^+}, (X^-)^{\prime,(3)}_{r^-}) \\
\quad\quad\quad\quad\quad\quad\quad\quad
\cdot(X^+)^{\prime,(2)}_{r^+}(X^-)^{\prime,(2)}_{r^-}, \\
X^+X^- =\sum_{r^+,r^-}
\tvt((X^+)^{\prime,(3)}_{r^+}, \tHopfS((X^-)^{\prime,(3)}_{r^-}))
\tvt((X^+)^{\prime,(1)}_{r^+}, (X^-)^{\prime,(1)}_{r^-}) \\
\quad\quad\quad\quad\quad\quad\quad\quad
\cdot(X^-)^{\prime,(2)}_{r^-}(X^+)^{\prime,(2)}_{r^+}
\end{array}
\end{equation} for $\lambda$, $\mu\in\defVZ$, $i$, $j\in\defI$,
and $X^+$, $Y^+\in\tBp$, $X^-$, $Y^-\in\tBm$,
where $(X^+)^{(x)}_{k^+}$ and $(X^-)^{(x)}_{k^-}$ with $x\in\bJ_{1,2}$
(resp. $(X^+)^{\prime,(y)}_{r^+}$ and
$(X^-)^{\prime,(y)}_{r^-}$ with $y\in\bJ_{1,3}$)
are any elements of $\tBp$ and $\tBm$ respectively
satisfying $\tHopfD(X^\pm)=\sum_{k^\pm}(X^\pm)^{(1)}_{k^\pm}
\otimes (X^\pm)^{(2)}_{k^\pm}$,
(resp. $((\rmid_\defU\otimes\tHopfD)\circ\tHopfD)(X^\pm)=
\sum_{r^\pm}(X^\pm)^{\prime,(1)}_{r^\pm}\otimes
(X^\pm)^{\prime,(2)}_{r^\pm}\otimes (X^\pm)^{\prime,(3)}_{r^\pm}$).

We have
\begin{equation*}
\begin{array}{l}
\tvtbhmPi(X^+\defacK_\lambda,X^-\defacL_\mu)=\tvtbhmPi(X^+,X^-)\bhm(\lambda,\mu)
\\
\quad (X^+\in\defUp,\,X^-\in\defUm,\,\lambda,\mu\in\defVZ).
\end{array}
\end{equation*}
It follows that $\tvtbhmPi_{|\defUp\times\defUm}$ is non-degenerate.
We have $\tvtbhmPi(X^+,X^-)=0$ for $\lambda,\mu\in\defVZPip$ with
$\lambda\ne\mu$ and $X^+\in\defUp_\lambda$,
$X^-\in\defUm_{-\mu}$.

\section{Kostant-Lusztig $\bA$-form}  \label{section:LAform}

In this section we establish the Kostant-Lusztig $\mathbb A$-forms for the 
finite-type multiparameter
quantum group $\defU=\defUbhmPi$ (see Definitions~\ref{definition:defofGQG} and \ref{definition:admbhm}) and construct invariant bases for each factor of the triangular decomposition. Let $A=[a_{ij}]_{ij\in\defI}$ be a Cartan matrix associated with a finite dimensional complex simple Lie algebra,
i.e., $A$ is a symmetrizable generalized Cartan matrix of finite-type. We assume that $\bhm$ is a $(\defPi,A)$-admissible bicharacter
(see Definition~\ref{definition:admbhm}). Recall the symbols $\tq_{ij}$, $\dotq_{ij}$
and $\defTheta(\tq_{ii}-1)$. Let $W$ be the Weyl group associated to the Cartan matrix $A$ and generated by the simple reflections
$\ts_i$ ($i\in\defI$) with relations
$\ts_i^2=\te$ ($i\in\defI$) and $(\ts_i\ts_j)^{\mbhmpi_{ij}}=\te$,
($\te$ is the identity element of $\tW$). Note that $W$ is the finite Weyl group.

\subsection{Some standard notations and results} \label{subsection:somno}
Define the map $\defell:W\to\bZp$ by $\defell(\te):=0$
and $\defell(w):=\min\{r\in\bN|\exists i_t\in\defI\,(t\in\bJ_{1,r}), w=\ts_{i_1}\cdots\ts_{i_r}\}$.
In fact $\defell$ is the length map of the Coxeter system $(\tW,\{\ts_i|i\in\defI\})$.
It is well-known that
\begin{equation*}
\mbox{$|\defell(w\ts_i)-\defell(w)|=1$ for $w\in\tW$ and $i\in\defI$.}
\end{equation*} Let $W$ act on $\defVZ$ by
$\ts_i\cdot\defPi(j):=\defPi(j)-a_{ij}\defPi(j)$ ($i,j\in\defI$). Use the convention as follows
\begin{equation*}
\mbox{Let $\ts_{i_1}\cdots\ts_{i_t}$ (resp. $\deftaubhm_{i_{t-1}}\cdots\deftaubhm_{i_1}\defPi$)
mean $\te$ (resp. $\defPi$) if $t=0$.}
\end{equation*}
We have
\begin{equation*}
\ts_{i_1}\cdots\ts_{i_{t-1}}\cdot\defPi(i_t)=\deftaubhm_{i_{t-1}}\cdots\deftaubhm_{i_1}\defPi(i_t)
\,\,(t\in\bN,\,i_x\in\defI\,(x\in\bJ_{1,t})).
\end{equation*}

For $w=\ts_{k_1}\ldots\ts_{k_r}\in\tW$, if $r=\defell(w)$,
the expression $\ts_{k_1}\ldots\ts_{k_{\defell(w)}}$ is called {\it{reduced}}.
For $w=\ts_{k_1}\ldots\ts_{k_{\defell(w)}}\in\tW$,
let
\begin{equation*}
\rmidbhmPiT_w:=\LusT_{k_1}^{\bhm,\defPi_1}\cdots\LusT_{k_{\defell(w)}}^{\bhm,\defPi_{\defell(w)}}\quad(\rmidbhmPiT_\te:=\rmid_{\defUbhmPi}),
\end{equation*}
where $\defPi_t:=\deftaubhm_{k_t}\cdots\deftaubhm_{k_1}\defPi$
($t\in\bJ_{1,{\defell(w)}}$). By Theorem~\ref{theorem:fundaTh}, $\rmidbhmPiT_w$ is independent of the
choice of reduced expressions for $w$.

 It is well-known that there exists a unique
$\twoh\in\tW$ such that $\defell(w)\leq\ellwo$ for all $w\in\tW$;
$\twoh$ is called {\it{the longest element}}.
We also know that $\ellwo=|\defRbhmpip|$ and
\begin{equation}\label{eqn:prtwoh}
\forall w\in\tW,\,\,\ellwo=\defell(w)+\defell(w^{-1}\twoh).
\end{equation}

Let $n=(n_1,\ldots,n_{\ellwo})\in\defI^{\ellwo}$ be such that $\ts_{n_1}\cdots\ts_{n_{\ellwo}}=\twoh$
(reduced expression of $\twoh$).
For $t\in\bJ_{1,\ellwo}$, let $\beta_{n;t}:=\ts_{n_1}\cdots\ts_{n_{t-1}}\cdot\defPi(n_t)$.
By
\eqref{eqn:HecIsoEq} and Lemma~\ref{lemma:lmfundaTh}, we have (see \cite[1.7]{Hum90})
\begin{equation*}
\defRbhmpip=\{\beta_{n;t}|t\in\bJ_{1,\ellwo}\}.
\end{equation*}
Thus $\defRbhmpi$ can be identified with the root system of $W$.
For $t\in\bJ_{1,\ellwo}$, let
\begin{equation*}
\begin{array}{l}
\defacE_{n;t}:=\rmidbhmPiT_{s_{n_1}\ldots s_{n_{t-1}}}(\defacE_{n_t}),\,\,
\defacF_{n;t}:=\rmidbhmPiT_{s_{n_1}\ldots s_{n_{t-1}}}(\defacF_{n_t}), \\
\BdefacE_{n;t}:=\rmidbhmPiT_{s_{n_1}\ldots s_{n_{t-1}}}(\BdefacE_{n_t}),\,\,
\BdefacF_{n;t}:=\rmidbhmPiT_{s_{n_1}\ldots s_{n_{t-1}}}(\BdefacF_{n_t}),
\end{array}
\end{equation*} ($\defacE_{n;1}:=\defacE_{n_1}$,
$\defacF_{n;1}:=\defacF_{n_1}$, $\BdefacE_{n;1}:=\BdefacE_{n_1}$,
$\BdefacF_{n;1}:=\BdefacF_{n_1}$), that is,
\begin{equation*}
\BdefacE_{n;t}={\frac {\defacE_{n;t}} {\defTheta(\tq_{n_tn_t}-1)}}
={\frac {\defacE_{n;t}} {\defTheta(\bhm(\beta_{n;t},\beta_{n;t})-1)}}.
\end{equation*}
Let $n_0\in\defI$ be such that $s_{n_0}\twoh s_{n_{\ellwo}}=\twoh$.
By Theorem~\ref{theorem:fundaTh}, using a standard argument (see \cite[Proposition~8.20]{Jan96}), we have
\begin{equation} \label{eqn:Enlwo}
\defacE_{n;\ellwo}=\defacE_{n_0},\,
\defacF_{n;\ellwo}=\defacF_{n_0},\,
\BdefacE_{n;\ellwo}=\BdefacE_{n_0},\,
\BdefacF_{n;\ellwo}=\BdefacF_{n_0}.
\end{equation}

The following result can be proved by  a standard argument (see \cite{AY2015}, \cite{HY10} for example).

\begin{theorem}  \label{theorem:LusPBW}
Let $k:=\ellwo$, $J:=\bJ_{1,k}$
and $\beta_t:=\ts_{n_1}\cdots\ts_{n_{t-1}}\cdot\defPi(n_t)$
$(t\in J)$.
\newline\newline
{\rm{(1)}} Let $\sigma:J\to J$
be a bijection. Then the elements
\begin{equation*}
\defacE_{n;\sigma(1)}^{x_1}\cdots\defacE_{n;\sigma(k)}^{x_k}\quad(x_t\in\bZp\,(t\in J)))
\end{equation*} form a $\bK$-basis of $\defUp$.
\newline
{\rm{(2)}} Let $y$, $z\in J$ be such that $y<z$.
Let $X$ be a $\bK$-subalgebra of $\defU$
generated by the elements $\defacE_{n;x}$ {\rm{(}}$x\in\bJ_{y+1,z-1}${\rm{)}}. Then
\begin{equation}\label{eqn:EnyEnz}
\defacE_{n;y}\defacE_{n;z}-\bhm(\beta_y,\beta_z)\defacE_{n;z}\defacE_{n;y}\in X.
\end{equation}
\newline
{\rm{(3)}} We have
\begin{equation*}
\tvtbhmPi(\defacE_{n;k}^{x_k}\cdots\defacE_{n;1}^{x_1},\defacF_{n;k}^{y_k}\cdots\defacF_{n;1}^{y_1})
=\prod_{t=1}^k\delta_{x_t,y_t}(x_t)_{\bhm(\beta_t,\beta_t)}!,
\end{equation*} where $x_t$, $y_t\in\bZp$ $(t\in J)$.
\end{theorem}

\subsection{Type ${\mathrm{G}}_2$} \label{subsection:typeGtwo}
In this subsection we assume $\defdimofV=2$ (so $\defI=\bJ_{1,2}$), and $A=[a_{ij}]_{ij\in\defI}$ is the Cartan
matrix of type $G_2$. So
 $a_{12}=-3$ and $a_{21}=-1$.
Consider the generalized quantum group $\defU=\defUbhmPi$.
Let $\tq_{ij}$ be as in Definition~\ref{definition:defofGQG} and $\bhm$ be a $(\defPi,A)$-admissible bicharacter.
Let $q:=\tq_{11}$ and $a:=\tq_{12}$.
Then $\tq_{22}=q^3$ and $\tq_{21}=q^{-3}a^{-1}$.
Since $\defacE_{4,1,2}=\defacE_{2,2,1}=0$, we have
\begin{equation*}
\begin{array}{l}
\defacE_1^4\defacE_2-(1+q)(1+q^2)a\defacE_1^3\defacE_2\defacE_1+q(1+q^2)(1+q+q^2)a^2\defacE_1^2\defacE_2\defacE_1^2 \\
\quad -q^3(1+q)(1+q^2)a^3\defacE_1^2\defacE_2\defacE_1^2+q^6a^4\defacE_2\defacE_1^4=0, \\
\defacE_1^2\defacE_2-(1+q)(1-q+q^2)a\defacE_2\defacE_1\defacE_2+q^3a^2\defacE_2^2\defacE_1=0.
\end{array}
\end{equation*}

Let
\begin{equation*}
\begin{array}{l}
\defacK_1;=\defacK_{\defPi(1)},\,\,\defacK_2;=\defacK_{\defPi(2)}, \,\,
 \defacE_{12}:=\defacE_{1,1,2},\,\,
 \defacE_{112}:=\defacE_{2,1,2},\\
 \defacE_{1112}:=\defacE_{3,1,2},\,\,
 \defacE_{11212}:=\defacE_{112}\defacE_{12}-aq^2\defacE_{12}\defacE_{112}.
\end{array}
\end{equation*}
Then we have
\begin{equation*}
\begin{array}{l}
\defacE_{12}\defacE_2=aq^3\defacE_2\defacE_{12},\,\,
\defacE_{1112}\defacE_2=a^3q^6\defacE_2\defacE_{1112}+a^2q^3(q^2-1)(q-1)\defacE_{12}^3, \\
\defacE_{112}\defacE_2=a^2q^3\defacE_2\defacE_{112}+aq(q^2-1)\defacE_{12}^2, \\
\defacE_{1112}\defacE_2=a^3q^3\defacE_2\defacE_{1112}+aq(q^2-q-1)\defacE_{11212}
+a^2q^2(q^3-1)\defacE_{12}\defacE_{112}, \\
\defacE_1\defacE_2=a\defacE_2\defacE_1+\defacE_{12},\,\,
\defacE_{11212}\defacE_{12}=aq^3\defacE_{12}\defacE_{11212}, \\
\defacE_{112}\defacE_{12}=aq^2\defacE_{12}\defacE_{112}+\defacE_{11212},\,\,
\defacE_{1112}\defacE_{12}=a^2q^3\defacE_{12}\defacE_{1112}+
{\frac {aq(q^3-1)} {q+1}}\defacE_{112}^2, \\
\defacE_1\defacE_{12}=aq\defacE_{12}\defacE_1+\defacE_{112},\,\,
\defacE_{112}\defacE_{11212}=aq^3\defacE_{112}\defacE_{11212}, \\
\defacE_{1112}\defacE_{11212}=a^3q^6\defacE_{11212}\defacE_{1112}+
{\frac {a^2q^3(q^3-1)(q-1)} {q+1}}\defacE_{112}^3, \\
\defacE_1\defacE_{11212}=a^2q^3\defacE_{11212}\defacE_1+
{\frac {aq(q^3-1)} {q+1}}\defacE_{112}^2,\,\,
\defacE_{1112}\defacE_{112}=aq^3\defacE_{112}\defacE_{1112}, \\
\defacE_1\defacE_{112}=aq^2\defacE_{112}\defacE_1+\defacE_{1112},\,\,
\defacE_1\defacE_{1112}=aq^3\defacE_{1112}\defacE_1.
\end{array}
\end{equation*}
We have
\begin{equation*}
\begin{array}{l}
\tHopfD(\defacE_{12})=\defacE_{12}\otimes 1+(1-q^{-3})\defacE_1\defacK_2\otimes\defacE_2
+\defacK_1\defacK_2\otimes\defacE_{12}, \\
\tHopfD(\defacE_{112})=\defacE_{112}\otimes 1+(1-q^{-3})(1-q^{-2})\defacE_1^3\defacK_2\otimes\defacE_2 \\
\quad\quad\quad\quad\quad +(1-q^{-2})(1+q)\defacE_1\defacK_1\defacK_2\otimes\defacE_{12}+\defacK_1^2\defacK_2\otimes\defacE_{112}, \\
\tHopfD(\defacE_{1112})=\defacE_{1112}\otimes 1+(1-q^{-3})(1-q^{-2})(1-q^{-1})\defacE_1^3\defacK_2\otimes\defacE_2 \\
\quad\quad\quad\quad\quad +(q^2-1)(1-q^{-3})\defacE_1^2\defacK_1\defacK_2\otimes\defacE_{12}, \\
\quad\quad\quad\quad\quad +q^{-1}(q^3-1)\defacE_1\defacK_1^2\defacK_2\otimes\defacE_{112}+\defacK_1^3\defacK_2\otimes\defacE_{1112},
\end{array}
\end{equation*}
and
\begin{equation*}
\begin{array}{l}
\tHopfD(\defacE_{11212})=\defacE_{11212}\otimes 1+{\frac {(q^3-1)^2} {q^4}}\defacE_{112}\defacE_1\defacK_2\otimes\defacE_2
+{\frac {(q^3-1)(q^2-q-1)} {aq^5}}\defacE_{1112}\defacK_2\otimes\defacE_2 \\
\quad\quad\quad\quad\quad\quad +{\frac {q^3-1} {q}}\defacE_{112}\defacK_1\defacK_2\otimes\defacE_{12}
+{\frac {(q^3-1)^2(q^2-1)(q-1)} {aq^{12}}}\defacE_1^3\defacK_1^2\otimes\defacE_2^2 \\
\quad\quad\quad\quad\quad\quad +{\frac {(q^3-1)^2(q^2-1)} {q^6}}\defacE_1^2\defacK_1\defacK_2\otimes\defacE_2\defacE_{12}
+{\frac {(q^3-1)(q^2-1)} {q^3}}\defacE_1\defacK_1^2\defacK_2^2\otimes\defacE_{12}^2 \\
\quad\quad\quad\quad\quad\quad +\defacK_1^3\defacK_2^2\otimes\defacE_{11212}.
\end{array}
\end{equation*}
Let
\begin{equation*}
\begin{array}{l}
\HdefacE_1:=\defacE_1,\,\HdefacE_2:=\defacE_2,\,\HdefacE_{12}:={\frac {q^3} {q^3-1}}\defacE_{12},
\,\HdefacE_{112}:={\frac {q^5} {(q^3-1)(q^2-1)}}\defacE_{112}, \\
\HdefacE_{1112}:={\frac {q^6} {(q^3-1)(q^2-1)(q-1)}}\defacE_{1112},\,
\HdefacE_{11212}:={\frac {q^9} {(q^3-1)^2(q^2-1)(q-1)}}\defacE_{11212}.
\end{array}
\end{equation*} Then we have
\begin{equation*}
\begin{array}{l}
\HdefacE_{1112}\HdefacE_2=a^3q^6\HdefacE_2\HdefacE_{1112}+a^2q^3(q^3-1)\HdefacE_{12}^3,\,
\HdefacE_{112}\HdefacE_2=a^2q^3\HdefacE_2\HdefacE_{112}+a(q^3-1)\HdefacE_{12}^2, \\
\HdefacE_{1112}\HdefacE_2=a^3q^3\HdefacE_2\HdefacE_{1112}+aq^{-2}(q^2-q-1)(q^3-1)\HdefacE_{11212} \\
\quad\quad\quad\quad\quad +a^2(q^3-1)(q^2+q+1)\HdefacE_{12}\HdefacE_{112}, \\
\HdefacE_1\HdefacE_2=a\HdefacE_2\HdefacE_1+q^{-3}(q^3-1)\HdefacE_{12},\,
\HdefacE_{112}\HdefacE_{12}=aq^2\HdefacE_{12}\HdefacE_{112}+q^{-1}(q-1)\HdefacE_{11212}, \\
\HdefacE_{1112}\defacE_{12}=a^2q^3\HdefacE_{12}\HdefacE_{1112}+a(q^3-1)\HdefacE_{112}^2,\,
\HdefacE_1\defacE_{12}=aq\HdefacE_{12}\HdefacE_1+q^{-2}(q^2-1)\HdefacE_{112}, \\
\HdefacE_{1112}\HdefacE_{11212}=a^3q^6\HdefacE_{11212}\HdefacE_{1112}+
a^2q^3(q^3-1)\HdefacE_{112}^3, \\
\HdefacE_1\HdefacE_{11212}=a^2q^3\HdefacE_{11212}\HdefacE_1+
a(q^3-1)\HdefacE_{112}^2,\,
\HdefacE_1\HdefacE_{112}=aq^2\HdefacE_{112}\HdefacE_1+q^{-1}(q-1)\HdefacE_{1112}.
\end{array}
\end{equation*} and
\begin{equation*}
\begin{array}{l}
\tHopfD(\HdefacE_{12})=\HdefacE_{12}\otimes 1+\HdefacE_1\defacK_2\otimes\HdefacE_2
+\defacK_1\defacK_2\otimes\HdefacE_{12}, \\
\tHopfD(\HdefacE_{112})=\HdefacE_{112}\otimes 1+\HdefacE_1^3\defacK_2\otimes\HdefacE_2
+(q+1)\HdefacE_1\defacK_1\defacK_2\otimes\HdefacE_{12}+\defacK_1^2\defacK_2\otimes\defacE_{112}, \\
\tHopfD(\HdefacE_{1112})=\HdefacE_{1112}\otimes 1+\HdefacE_1^3\defacK_2\otimes\HdefacE_2
+(q^2+q+1)\HdefacE_1^2\defacK_1\defacK_2\otimes\HdefacE_{12}, \\
\quad\quad\quad\quad\quad +(q^2+q+1)\HdefacE_1\defacK_1^2\defacK_2\otimes\HdefacE_{112}+\defacK_1^3\defacK_2\otimes\HdefacE_{1112}, \\
\tHopfD(\HdefacE_{11212})=\HdefacE_{11212}\otimes 1
+(q^2+q+1)\HdefacE_{112}\HdefacE_1\defacK_2\otimes\HdefacE_2 \\
\quad\quad\quad\quad\quad\quad +a^{-1}q^{-2}(q^2-q-1)\HdefacE_{1112}\defacK_2\otimes\HdefacE_2 \\
\quad\quad\quad\quad\quad\quad +(q^2+q+1)\HdefacE_{112}\defacK_1\defacK_2\otimes\HdefacE_{12}
+a^{-1}q^{-3}\HdefacE_1^3\defacK_1^2\otimes\HdefacE_2^2 \\
\quad\quad\quad\quad\quad\quad +(q^2+q+1)\HdefacE_1^2\defacK_1\defacK_2\otimes\HdefacE_2\HdefacE_{12}
+(q^2+q+1)\HdefacE_1\defacK_1^2\defacK_2^2\otimes\HdefacE_{12}^2 \\
\quad\quad\quad\quad\quad\quad +\defacK_1^3\defacK_2^2\otimes\HdefacE_{11212}.
\end{array}
\end{equation*}

For ${\bf a}=(a_1,\ldots,a_6)\in\bZp^6$, let
\begin{equation*}
\begin{array}{l}
\defQp({\bf a}):=
{\frac {\defacE_2^{a_1}
\defacE_{12}^{a_2}
((3)_q^{-1}!\defacE_{11212})^{a_3}
((2)_q^{-1}\defacE_{112})^{a_4}
((3)_q^{-1}!\defacE_{1112})^{a_5}
\defacE_1^{a_6}}
{(a_1)_{q^3}!(a_2)_q!(a_3)_{q^3}!(a_4)_q!(a_5)_{q^3}!(a_6)_q!}
}, \\
\defQm({\bf a}):=\HdefacE_2^{a_1}\HdefacE_{12}^{a_2}\HdefacE_{11212}^{a_3}\HdefacE_{112}^{a_4}
\HdefacE_{1112}^{a_5}\HdefacE_1^{a_6}.
\end{array}
\end{equation*}
Then we have
\begin{equation}\label{eqn:dualGtwo}
\begin{array}{l}
\tvtbhmPi(\defQp({\bf a}),\defUpbhmopPi(\defQm({\bf b})))
=\delta_{{\bf a},{\bf b}}\quad({\bf a},\,{\bf b}\in\bZp^6).
\end{array}
\end{equation}
By \eqref{eqn:dualGtwo},
$\{\defQp({\bf a})|{\bf a}\in\bZp^6\}$ and $\{\defQm({\bf a})|{\bf a}\in\bZp^6\}$
are $\bK$-bases of $\defUbhmPi$.
Let ${\acute{\bA}}$ be the $\bZ$-subalgebra of $\bK$ generated by $q^{\pm 1}$
and $a^{\pm 1}$, i.e., ${\acute{\bA}}=\bZ[q^{\pm 1},a^{\pm 1}]$.
For $t\in\bJ_{1,2}$, let $\mfkgpm$ be the ${\acute{\bA}}$-submodules
of $\tBp$ with the ${\acute{\bA}}$-bases
$\{\defQpm({\bf a})\defacK_\lambda|{\bf a}\in\bZp^6,\,\lambda\in\defVZ\}$.
Clearly $\mfkgm$ is a ${\acute{\bA}}$-subalgebra of $\defU$
with $\tHopfD(\mfkgm)\subset\mfkgm\otimes_{{\acute{\bA}}}\mfkgm$.
By \eqref{eqn:bprtoftB} and \eqref{eqn:dualGtwo}, we see the following.
\begin{equation}\label{eqn:Key}
\begin{array}{l}
\mbox{As a ${\acute{\bA}}$-subalgebra of $\defU$, $\mfkgp$ is generated by
${\frac {\defacE_i^x} {(x)_{\tq_{ii}}!}}$} \\
\mbox{($i\in\defI(=\bJ_{1,2})$, $x\in\bZp$) and $\defacK_\lambda$ ($\lambda\in\defVZ$).} \\
\mbox{$\mfkgp$ is a Hopf ${\acute{\bA}}$-subalgebra of $\defU$.}
\end{array}
\end{equation} Then $\mfkgm$ is also a Hopf ${\acute{\bA}}$-algebra.

Let $n:=(1,2,1,2,1,2)\in\defI^{\ellwo}$.
Then we have
\begin{equation}\label{eqn:Keyeqa}
\begin{array}{ll}
\BdefacE_{n;1}=\BdefacE_1, &  \BdefacF_{n;1}=\BdefacF_1(=\defUpbhmopPi(\BdefacE_{n;1})), \\
\BdefacE_{n;2}={\frac{\dotq_{12}^{-3}} {(3)_q!\defTheta(q-1)^3\defTheta(q^3-1)}}\defacE_{1112},    &
\BdefacF_{n;2}=\dotq_{11}^3\defUpbhmopPi(\BdefacE_{n;2}), \\
\BdefacE_{n;3}={\frac{\dotq_{12}^{-2}} {(2)_q!\defTheta(q^3-1)\defTheta(q-1)^2}}\defacE_{112}, &
\BdefacF_{n;3}=\dotq_{11}^4\defUpbhmopPi(\BdefacE_{n;3}), \\
\BdefacE_{n;4}={\frac{\dotq_{12}^{-4}} {(3)_q!\defTheta(q^3-1)^2\defTheta(q-1)^3}}\defacE_{11212}, &
\BdefacF_{n;4}=\dotq_{11}^6\defUpbhmopPi(\BdefacE_{n;4}), \\
\BdefacE_{n;5}={\frac{\dotq_{12}^{-1}} {\defTheta(q-1)\defTheta(q^3-1)}}\defacE_{12}, &
\BdefacF_{n;5}=\dotq_{11}^3\defUpbhmopPi(\BdefacE_{n;5}), \\
\BdefacE_{n;6}=\BdefacE_2, & \BdefacF_{n;6}=\BdefacF_2(=\defUpbhmopPi(\BdefacE_{n;6})),
\end{array}
\end{equation} where recall that $\dotq_{ij}^2=a$.
Let $n^\prime:=(2,1,2,1,2,1)\in\defI^{\ellwo}$. By \eqref{eqn:InvLusTi} and  \eqref{eqn:Enlwo}, we have
\begin{equation}\label{eqn:Keyeqb}
\BdefacE_{n^\prime;t}=\GammabhmopPi(\BdefacE_{n;7-t}),\,\,\BdefacF_{n^\prime;t}=\GammabhmopPi(\BdefacF_{n;7-t})
\quad (t\in\bJ_{1,6}).
\end{equation}

Let $\barmfkgp$ be the ${\acute{\bA}}$-subalgebra of $\defU$ generated by
${\frac {\BdefacE_i^x} {(x)_{\tq_{ii}}!}}$ ($i\in\defI(=\bJ_{1,2})$, $x\in\bZp$) and $\defacK_\lambda$ ($\lambda\in\defVZ$).
Let $y_{2t-1}:=\tq$, $y_{2t}:=\tq^3$ for $t\in\bJ_{1,3}$.
Let $y^\prime_t:=y_{7-t}$  for $t\in\bJ_{1,6}$.
By \eqref{eqn:EnyEnz}, \eqref{eqn:Key}, \eqref{eqn:Keyeqa} and \eqref{eqn:Keyeqb}, we see the following.
\begin{equation}\label{eqn:BdefEKey}
\begin{array}{l}
\mbox{For a bijection $\sigma:\bJ_{1,6}\to\bJ_{1,6}$, the elements} \\
\mbox{${\frac {\BdefacE_{n;\sigma(1)}^{x_1}} {(x_1)_{y_{\sigma(1)}}!}}\cdots{\frac {\BdefacE_{n;\sigma(6)}^{x_6}} {(x_6)_{y_{\sigma(6)}}!}}\defacK_\lambda$
(resp. ${\frac {\BdefacE_{n^\prime;\sigma(1)}^{x_1}} {(x_1)_{y^\prime_{\sigma(1)}}!}}\cdots{\frac {\BdefacE_{n^\prime;\sigma(6)}^{x_6}} {(x_6)_{y^\prime_{\sigma(6)}}!}}\defacK_\lambda$)} \\
\mbox{($x_t\in\bZp$ ($t\in\bJ_{1,6}$), $\lambda\in\defVZ$) form a ${\acute{\bA}}$-base of $\barmfkgp$.} \\
\mbox{$\barmfkgp$ is a Hopf ${\acute{\bA}}$-subalgebra of $\defU$.}
\end{array}
\end{equation}

Let $\barmfkgm$ be the ${\acute{\bA}}$-subalgebra of $\defU$ generated by
${\frac {\BdefacF_i^x} {(x)_{\tq_{ii}}!}}$ ($i\in\defI(=\bJ_{1,2})$, $x\in\bZp$) and $\defacL_\lambda$ ($\lambda\in\defVZ$).
By \eqref{eqn:Keyeqa}, \eqref{eqn:Keyeqb}
and \eqref{eqn:BdefEKey}, we see the following.
\begin{equation}\label{eqn:BdefFKey}
\begin{array}{l}
\mbox{For a bijection $\sigma:\bJ_{1,6}\to\bJ_{1,6}$, the elements} \\
\mbox{${\frac {\BdefacF_{n;\sigma(1)}^{x_1}} {(x_1)_{y_{\sigma(1)}}!}}\cdots{\frac {\BdefacF_{n;\sigma(6)}^{x_6}} {(x_6)_{y_{\sigma(6)}}!}}\defacL_\lambda$
(resp. ${\frac {\BdefacF_{n^\prime;\sigma(1)}^{x_1}} {(x_1)_{y^\prime_{\sigma(1)}}!}}\cdots{\frac {\BdefacF_{n^\prime;\sigma(6)}^{x_6}} {(x_6)_{y^\prime_{\sigma(6)}}!}}\defacL_\lambda$)} \\
\mbox{($x_t\in\bZp$ ($t\in\bJ_{1,6}$), $\lambda\in\defVZ$) form a ${\acute{\bA}}$-base of $\barmfkgm$.} \\
\mbox{$\barmfkgm$ is a Hopf ${\acute{\bA}}$-subalgebra of $\defU$.}
\end{array}
\end{equation}

\subsection{General case}
Let $\bA$ be the $\bZ$-subalgebra of $\bK$ generated by
$\dotq_{ij}^{\pm 1}$ for all $i$, $j\in\defI$.
$n=(n_1,\ldots,n_{\ellwo})\in\defI^{\ellwo}$ be as in Subsection~\ref{subsection:somno}.
Let $\defPi_0:=\defPi$ and
$\defPi_t:=\deftaubhm_{n_t}\cdots\deftaubhm_{n_1}\defPi$
($t\in\bJ_{1,\ellwo}$).
For $t\in\bJ_{1,\kappa}$ and $x\in\bZp$, let
\begin{equation*}
\BdefacE_{n;t}^{(x)}:={\frac {\BdefacE_{n;t}^x}
{(x)_{\bhm(\defPi_{t-1}(n_t),\defPi_{t-1}(n_t))}!}},\quad
\BdefacF_{n;t}^{(x)}:={\frac {\BdefacF_{n;t}^x}
{(x)_{\bhm(\defPi_{t-1}(n_t),\defPi_{t-1}(n_t))}!}}.
\end{equation*}
Let $\defUpA$ (resp. $\defUmA$) be the $\bA$-subalgebra of $\defUp$
(resp. $\defUm$) generated by
${\frac {\BdefacE_i^x} {(x)_{\tq_{ii}}!}}$
(resp. ${\frac {\BdefacF_i^x} {(x)_{\tq_{ii}}!}}$)
with $i\in\defI$ and $x\in\bZp$.
\begin{theorem}
Let $\sigma:\bJ_{1,\ellwo}\to\bJ_{1,\ellwo}$ be a bijection.
\newline\newline
{\rm{(1)}}
The elements
\begin{equation} \label{eqn:mainUpA}
\mbox{$\BdefacE_{n;\sigma(1)}^{(x_1)}\cdots\BdefacE_{n:\sigma(\ellwo)}^{(x_{\ellwo})}$
with $x_t\in\bZp$ {\rm{(}}$t\in\bJ_{1,\ellwo}${\rm{)}} }
\end{equation} form an
$\bA$-base for $\defUpA$.
\newline
{\rm{(2)}}
The elements
\begin{equation} \label{eqn:mainUmA}
\mbox{$\BdefacF_{n;\sigma(1)}^{(x_1)}\cdots\BdefacF_{n:\sigma(\ellwo)}^{(x_{\ellwo})}$
with $x_t\in\bZp$ {\rm{(}}$t\in\bJ_{1,\ellwo}${\rm{)}} }
\end{equation} form an
$\bA$-base for $\defUmA$.
\end{theorem}
{\it{Proof.}} We only prove (1), as (2) can be proved similarly or obtained from (1) via Chevalley isomorphism.
If $|I|=1$, the claim is clear.

If $|I|=2$ and $A$ is of type $G_2$, the result is obtained in \eqref{eqn:BdefEKey} and \eqref{eqn:BdefFKey}.
For other rank two types, the claim can be proved in a similar and in fact easier way.

Let us consider the higher rank cases.
Let $Z$ be the free $\bA$-submodule of $\defUp$ with the
free $\bA$-basis formed by the elements of \eqref{eqn:mainUpA} (for $\sigma:=\rmid_{\bJ_{1,\ellwo}}$,
see Theorem~\ref{theorem:LusPBW}~(1)).
By an argument similar to \cite[Proposition~8.20]{Jan96}
and by the claim for rank-two cases, we see that
$\BdefacE_{n;t}^{(x)}\in\defUpA$, so $Z\subseteq\defUpA$.
Now by \eqref{eqn:prtwoh}, we have $Z=\defUpA$.
Thus the claim follows from \eqref{eqn:EnyEnz}. \hfill $\Box$
\newline\par
For $k$, $r\in\bZp$ and $x\in\bKt$
with $(k+r)_x!\ne 0$, 
let ${{k} \choose {r}}_x:={\frac {(k+r)_x!} {(k)_x!(r)_x!}}$.
For $x\in\bKt$ with $x^b\ne 1$ for all $b\in\bN$ and for $X$, $Y\in\defUo$, $l\in\bZ$ and $p\in\bZp$, let
\begin{equation*}
\left[\begin{array}{c} X,Y,l \\ p\end{array}\right]_x
:=\prod_{t=1}^p{\frac {x^{l-t+1}X-Y} {x^t-1}}.
\end{equation*} Then we have
\begin{equation*}
\begin{array}{l}
\left[\begin{array}{c} X,Y,l \\ p\end{array}\right]_x
-\left[\begin{array}{c} X,Y,l+1 \\ p\end{array}\right]_x
=-x^{l-p+1}\left[\begin{array}{c} X,Y,l \\ p-1\end{array}\right]_xX, \\
\left[\begin{array}{c} X,Y,0 \\ l\end{array}\right]_x
\left[\begin{array}{c} X,Y,-l \\ p\end{array}\right]_x
={{p+l} \choose {p}}_x\left[\begin{array}{c} X,Y,0 \\ p+l\end{array}\right]_x.
\end{array}
\end{equation*}
Let $\defUoA$ be the $\bA$-subalgebra generated by
$\defacK_\lambda\defacL_\mu$ ($\lambda$, $\mu\in\defVZ$)
and $\left[\begin{array}{c} \defacK_{\defPi(i)},\defacL_{\defPi(i)},l \\ p\end{array}\right]_{\tq_{ii}}$
($i\in\defI$, $l\in\bZ$, $p\in\bZp$).
By a standard argument, we have the following two lemmas.
\begin{lemma} \label{lemma:bsUoA} The elements
\begin{equation*}
\prod_{i\in\defI}\defacK_{\defPi(i)}^{x_i}(\defacK_{\defPi(i)}\defacL_{\defPi(i)})^{y_i}
\left[\begin{array}{c} \defacK_{\defPi(i)},\defacL_{\defPi(i)},0 
\\ z_i\end{array}\right]_{\tq_{ii}}\quad
(x_i\in\bZ_{0,1},\,y_i\in\bZ,\,z_i\in\bZp)
\end{equation*} form an $\bA$-base of $\defUoA$.
\end{lemma}

Let $\defUA$ be the $\bA$-subalgebra of $\defU$ generated by
$\defUpA$ and $\defUmA$.

\begin{lemma} \label{lemma:triUA}
We have $\defUoA\subset\defUA$, and we have
the $\bA$-module isomorphism
\begin{equation*}
{\mathbf{m}}_\bA:\defUmA\otimes_\bA\defUoA\otimes_\bA\defUpA\longrightarrow \defUA
\end{equation*} given by ${\mathbf{m}}_\bA(X\otimes Y\otimes Z):=XYZ$.
\end{lemma}

\begin{center}
{\bf{Acknowledgments}}
\end{center}

We thank the referee for careful reading and kind comments. The third author would like to express his heartfelt thanks to Professor Naihong~Hu for valuable communication.

\end{document}